%% file: main.tex
\documentclass[preprint]{elsarticle} 

\makeatletter
\def\ps@pprintTitle{%
\let\@oddhead\@empty
\let\@evenhead\@empty
\def\@oddfoot{\centerline{\thepage}}%
\let\@evenfoot\@oddfoot}
\makeatother

\usepackage{etoolbox}
\patchcmd{\MaketitleBox}{\footnotesize\itshape\elsaddress\par\vskip36pt}{\footnotesize\itshape\elsaddress\par\parbox[b][36pt]{\linewidth}{\vfill\hfill\textnormal{\today}\hfill\null\vfill}}{}{}%
\patchcmd{\pprintMaketitle}{\footnotesize\itshape\elsaddress\par\vskip36pt}{\footnotesize\itshape\elsaddress\par\parbox[b][36pt]{\linewidth}{\vfill\hfill\textnormal{\today}\hfill\null\vfill}}{}{}%

\usepackage[
colorlinks=true,%
breaklinks=true,%
linkcolor=blue,
urlcolor=blue,%
citecolor=blue,%
pdftitle={Risk-based Design Optimization for Powder Bed Fusion Metal Additive Manufacturing}, 
pdfkeywords={},	
pdfauthor={Yulin Guo},		
bookmarksopen=false,
pdfpagemode=UseNone]{hyperref}

\usepackage[margin = 1.2in]{geometry}
\usepackage[T1]{fontenc}
\usepackage[english]{babel}
\usepackage[latin1]{inputenc}
\usepackage{mathtools}
\usepackage{amsmath}
\usepackage{amsfonts}
\usepackage{amsthm}
\usepackage{amssymb}
\usepackage{array}

\usepackage{subcaption}

\usepackage{siunitx}

\usepackage{color}
\usepackage{url}
\usepackage{cleveref}
\usepackage{graphicx}
\usepackage{breakcites}
\usepackage{enumitem}
\usepackage[title,titletoc,toc]{appendix}

\usepackage[textsize=tiny]{todonotes}

\usepackage{soul}
\usepackage{algorithm}
\usepackage{algorithmicx}
\usepackage{algpseudocode}
\usepackage{tabularx}
\usepackage{array}
\usepackage{multirow}
\usepackage{threeparttable}
\usepackage{booktabs}
\usepackage{graphicx}
\usepackage{arydshln}
\usepackage{setspace}   
\usepackage{xcolor}
\usepackage[mathlines,running]{lineno}
\usepackage{epstopdf}
\usepackage{tikz}
\usepackage{tikz-imagelabels}
\usepackage{pgfplots} 
\usepackage{pgfgantt}
\usepgfplotslibrary{groupplots}
\pgfplotsset{compat=1.17}

\input{mymacros.tex}

\newcommand\Algphase[1]{%
\vspace*{-.4\baselineskip}\Statex\hspace*{\dimexpr-\algorithmicindent-2pt\relax}\rule{\textwidth}{0.4pt}%
\Statex\hspace*{-\algorithmicindent}\textbf{#1}%
\vspace*{-.7\baselineskip}\Statex\hspace*{\dimexpr-\algorithmicindent-2pt\relax}\rule{\textwidth}{0.4pt}%
}

\begin{document}	

\begin{frontmatter}
		
\title{Risk-based Design Optimization for Powder Bed Fusion Metal Additive Manufacturing}

\author[ucsd]{{Yulin Guo}}
\ead{yug054@ucsd.edu}

		\author[ucsd]{Boris Kramer\corref{cor1}}
		\ead{bmkramer@ucsd.edu}

		\cortext[cor1]{Corresponding author}
						
		\address[ucsd]{Department of Mechanical and Aerospace Engineering, University of California San Diego, CA, United States}

\begin{abstract}
Powder bed fusion is a widely used additive manufacturing (AM) process for producing complex, small-batch parts that are impractical to manufacture using conventional methods. However, its broader adoption is hindered by process-induced defects. The challenge in AM stems from inherent material and process uncertainties. Therefore, it is critical to account for these uncertainties in the design optimization and control of powder bed fusion AM processes. In this work, we formulate and solve a design optimization problem under uncertainty for a powder bed fusion metal AM process. Our objective is to minimize energy consumption while enforcing a risk-based constraint formulated with a buffered probability of failure on residual stress, along with a constraint on melting temperature to ensure a successful build. We use surrogate models for the residual stress and temperature snapshots to accelerate optimization; we train these models using data from high-fidelity finite element simulations. We validate the optimization results through additional high-fidelity simulations. The validated results demonstrate that the proposed optimization reduces energy consumption, enhances process reliability, and contributes to more robust and sustainable additive manufacturing.
\end{abstract}

		
\end{frontmatter}
\section{Introduction}


Additive manufacturing (AM) is the process of joining materials to make parts from 3D model data, usually layer upon layer, as opposed to subtractive manufacturing and formative manufacturing methodologies~\cite{cdi_astm_standards_111828, gibson2021additive}. It has grown in popularity due to its speed, low labor cost, customization ability, and ability to make highly complex geometric designs compared to traditional manufacturing methodologies, such as milling and injection molding. A variety of materials, including polymers, metals, and ceramics, are available to provide the unique attributes needed to satisfy the engineering requirements of AM parts in aerospace, automotive, medical devices, and other industries~\cite{shapiro2016additive, vasco2021additive, seoane2021semi, murr2018additive, wang2024uncertainty}. The American Society for Testing and Materials classifies AM processes into seven categories: binder jetting, directed energy deposition, material extrusion, material jetting, powder bed fusion, sheet lamination, and vat photopolymerization~\cite{cdi_astm_standards_111828}. Several additive manufacturing processes are available for metallic materials, with each process utilizing specific material forms: metal sheets for sheet lamination, metal rods or wires for material extrusion, metal wires for directed energy deposition, and metal powders for binder jetting and powder bed fusion~\cite{konda2017additive}. Powder bed fusion methods have been applied to a wide range of metallic materials such as Fe-based alloys~\cite{wu2014experimental}, Ni-based alloys~\cite{lu2015study}, Ti-based alloys~\cite{gong2015influence}, Co-based alloys~\cite{wei2018heterogeneous}, and Cu-based alloys~\cite{lu2021additive}. Within the powder bed fusion methods, electron beam melting and selective laser melting (also known as laser beam melting) are two common methods to produce metallic parts where the powder is melted and fused by an electron or laser beam. These two processes share the same fundamental workflow with differences in their operating environments. In electron beam melting, a heated powder bed and a vacuum chamber are used. The vacuum guides and focuses the electron beam. By comparison, selective laser melting operates with a cold powder bed inside a closed chamber filled with inert gases such as N$_2$ or Ar, which prevents oxidation of the metal powder during melting~\cite{fu20143, konda2017additive}. 


The selective laser melting process involves over 50 parameters, including laser power, scanning speed, hatch distance, and overlaps; the electron beam melting process encompasses even more parameters, such as beam focus, beam diameter, beam line spacing, plate temperature, preheating temperature, contour strategy, and scan strategy~\cite{spears2016process, konda2017additive}.  Improper selection of processing parameters can lead to defects such as high surface roughness, pores and voids, cracks, delaminations, or distortions, which significantly compromises part quality~\cite{gong2015influence, yan2017multi, bartlett2019overview, galarraga2016effects}. For example, higher scanning speeds increase the length-to-diameter ratio of the molten pool present in the manufacturing process. When this ratio exceeds $\pi$, it causes the balling effect (breaking up of the molten pool into spherical particles). This phenomenon degrades the surface integrity of the final components. Scanning speed and cooling rate can also affect the grain size, which in turn affects the mechanical properties~\cite{corbin2018effect, yadroitsev2013energy}.

 
The relationships between process parameters and part quality involve complex nonlinear correlations that do not allow for explicit mathematical formulation~\cite{cao2021optimization}. Despite the growing adoption of AM technologies, the field lacks comprehensive guidelines for selecting appropriate AM processes based on material requirements and desired part properties~\cite{konda2017additive}. Process parameter optimization for AM processes has been extensively studied using design of experiments (DoE) methods, including factorial design~\cite{gong2014analysis, oyesola2021optimization}, central composite design~\cite{bhardwaj2019direct, elsayed2018optimization}, and Taguchi methods~\cite{sun2013parametric, gong2014analysis}. These approaches are often complemented by statistical analysis techniques such as regression analysis~\cite{gong2014analysis, gong2015influence, shi2017parameter} and analysis of variance (ANOVA)~\cite{sun2013parametric, elsayed2018optimization}. These optimization methods are time-consuming because they typically require evaluating multiple parameter combinations. The computational burden grows exponentially as each additional parameter requires testing across multiple intervals, either through numerical simulations or physical experiments. To address these limitations, AM process parameter optimization using surrogate models has been proposed in literature, with studies incorporating Gaussian process modeling~\cite{meng2020process, cao2021optimization} and various machine learning approaches~\cite{gaikwad2020heterogeneous}. 


The uncertainties in AM process parameters propagate through the manufacturing process, making quantities of interest, such as the residual stress of the printed part, uncertain as well. Consequently, design for AM processes should consider uncertainties not only in the process parameters but also in the quantities related to manufactured part qualities, such as residual stress. In order to incorporate uncertainties, robust design optimization~\cite{lee2021robust} and reliability-based design optimization~\cite{zou2006direct} have been used. The former aims to minimize the sensitivity of the system performance and typically uses mean and variance as optimization objectives. It does not explicitly account for failure probability, which may result in designs that are consistent but still have unacceptable failure rates in extreme conditions. A robust design optimization approach has been proposed for process parameters in metallic additive manufacturing~\cite{wang2019data}. Their work addresses uncertainties in material properties and power absorption coefficients. Reliability-based design optimization manages uncertainty by enforcing the probability of failure below specified thresholds while optimizing performance objectives. It uses a first- or second-order reliability method or Monte Carlo simulation to evaluate failure. A major limitation is the high computational cost, especially for estimating rare-event probabilities. A traditional approach is to add safety margins to the results of deterministic design optimization to compensate for uncertainties. However, this approach is inefficient, as it neither incorporates uncertainties directly into the design optimization process nor optimizes the safety margins themselves, which can potentially limit performance. Heuristic optimization methods have also been used in recent studies, such as whale optimization algorithm in~\cite{cao2021optimization}. A principled way to account for the failure magnitude integrates superquantile (also known as conditional value-at-risk, CVaR)~\cite{rockafellar2002conditional} and buffered probability of failure (BPOF)~\cite{rockafellar2010buffered} into design optimization processes. These alternative risk measures can be used either as objectives or constraints in optimization formulations~\cite{chaudhuri2022certifiable} and eliminate the guesswork associated with choosing safety margins. 


The main contribution of this article is to formulate and solve a risk-based design optimization problem for a powder bed fusion metal AM process under uncertainty, enabling simultaneous reduction in energy consumption and improvements in process reliability. This advances robust and sustainable metal additive manufacturing by explicitly quantifying and mitigating failure risks while optimizing process parameters. We achieve this as follows: (a) we incorporate a buffered probability of failure~(BPOF) constraint on the residual stress of the manufactured part to ensure the optimization penalizes both the frequency and severity of process failures. (b) We transform the temperature history and residual stress of the manufactured part into low-dimensional spaces, in which we train surrogate models with high-fidelity thermal-mechanical simulation data. This enables efficient optimizations. (c) We validate the optimization results through additional high-fidelity simulations and confirm reductions in energy consumption without compromising part quality.

The rest of this paper is organized as follows. In Section~\ref{section: problem_discription}, we present an overview of risk-based design optimization. In Section~\ref{section: background}, we introduce the powder bed fusion process, its thermal and mechanical governing equations, finite element modeling, residual stress development in manufactured parts, and mitigation strategies. In Section~\ref{section: DUU_AM}, we then present the design under uncertainty problem in full detail and cover risk-based measures, surrogate model development, and our optimization methodology. In section~\ref{section: numerical_results}, we present and analyze numerical results. Finally, we summarize our conclusions and suggest directions for future research in Section~\ref{section: conclusion}.
%
\section{Formulation of risk-based design optimization for the powder bed fusion metal additive manufacturing problem} \label{section: problem_discription} 

In powder bed fusion, the beam power and scanning speed have emerged as predominant optimization parameters across recent studies~\cite{meng2020process, shi2017parameter, cao2021optimization}. We formulate an optimization problem with the goal to minimize the total energy consumed, which is a function of the beam power and scanning speed, subject to a threshold on the BPOF on the maximum residual stress, as well as a threshold on the maximum temperature, in a part manufactured by an electron beam melting process. 

The manufacturing process is parameterized by two design variables, four random variables and six parameters. The two design variables are the scanning speed $v \in [v_L, v_U]$ and beam power $P \in [P_L, P_U]$, which form the design vector $\bd = [v, P]^{\top}$. The four random variables are the preheating temperature $T_0$ of the machine, the yield strength $Y$, elastic modulus $E$, and bulk density $\rho$ of the material. These variables are assumed to be statistically independent and are collected in the random vector $Z = [T_0, Y, E, \rho]^{\top}$. The specific heat $C_p(T)$ and thermal conductivity $\kappa(T)$ are functions of temperatures and are parameterized by coefficients, $a_{i}$ and $b_{i}, \ i = 1,2,3$, respectively, as discussed in detail in Section~\ref{section: bg_thermal_model} and Section~\ref{section: simulation_runs}. The six coefficients are taken to be deterministic parameters and are collected in a vector $\boldsymbol{\theta}=[a_1,a_2, a_3,b_1,b_2,b_3]^{\top}$, which we do not optimize over. We formulate the optimization problem as follows:
\begin{equation}
    \begin{aligned}
    & \underset{\substack{P_L \leq P \leq P_U, \\ v_L \leq v \leq v_U}}{\text{min}}
    & & P \cdot \frac{l}{v} \\
    & \text{subject to}
    & & \bar{p}(\sigma_{\max}(\bd, Z; \boldsymbol{\theta})) \leq 1 - \alpha_T,\\
    & & & T_{\text{liq}} < T_{\max}(\bd, Z; \boldsymbol{\theta}) < 1.1 \times T_{\text{liq}},
    \end{aligned}
    \label{eq: optimization_setup}
\end{equation} 
\noindent where $\bar{p}(\sigma_{\max}(\bd, Z;\boldsymbol{\theta}))$ is the BPOF, $\alpha_T$ is the desired reliability level, $\sigma_{\max}$ is the maximum residual stress, $T_{\max}$ is the maximum temperature in the temperature history, and $T_{\text{liq}}$ is the liquidus temperature of the material, see Section~\ref{section: simulation_runs} for more details and specific values. 

Next, we discuss the powder bed fusion application, including its governing equations to determine $\sigma_{\max}$ and $T_{\max}$. Then in Section~\ref{section: DUU_AM}, we discuss the reliability constraint on $\sigma_{\max}$ followed by the solution to the optimization problem~\eqref{eq: optimization_setup}.
%
\section{Powder bed fusion application} \label{section: background}
%
We introduce the electron beam melting process in Section~\ref{section: PBF} and the residual stress that arises in additively manufactured parts in Section~\ref{section: bg_residual_stress}. We then discuss the thermal and mechanical models governing the powder bed fusion process in Section~\ref{section: bg_thermal_model} and Section~\ref{section: bg_mechanical_model}. Finally, we present the finite element models used to simulate the manufacturing process in Section~\ref{section: FEM}. 

\subsection{Electron beam melting process} \label{section: PBF}

\begin{figure}[htbp]
    \centering
    \begin{subfigure}{0.9\textwidth}
        \includegraphics[width=\textwidth]{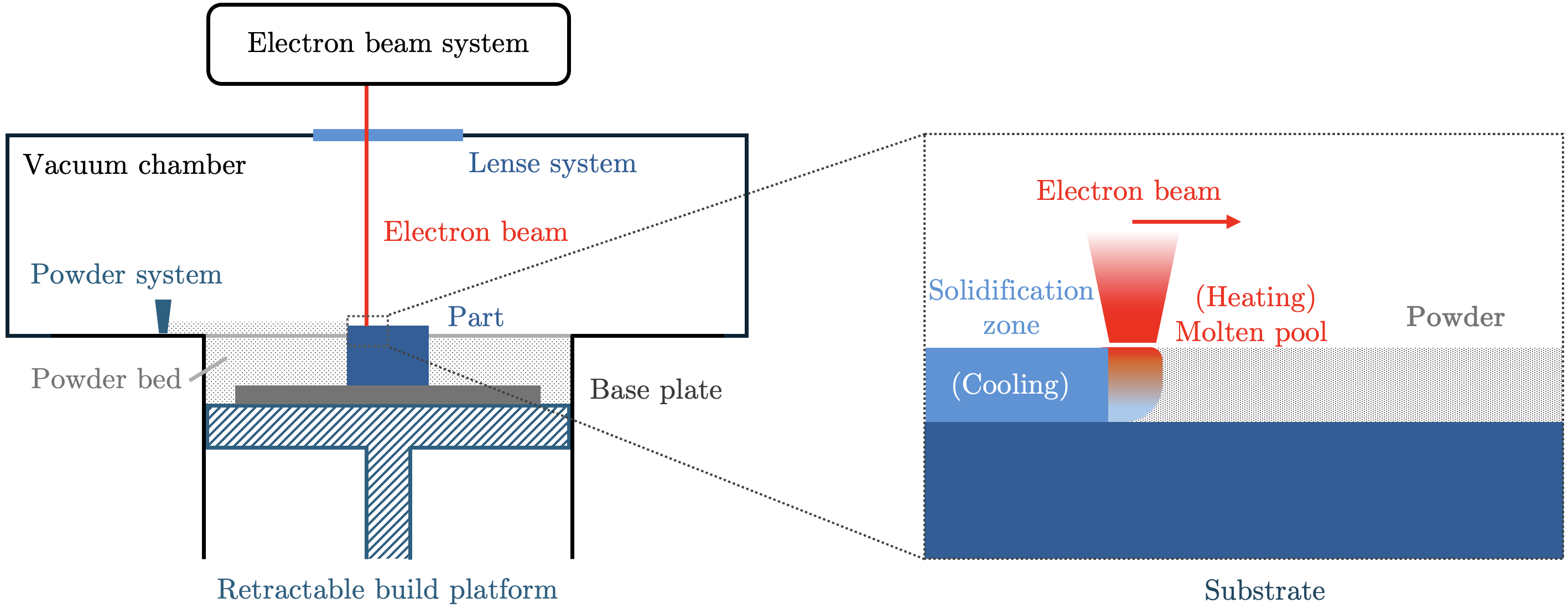}
    \end{subfigure}
    \caption{Schematic of an electron beam melting machine and the printing process.}
    \label{fig: EBM}
\end{figure}

In an electron beam melting powder bed fusion process, a thin layer of powder is deposited over a substrate plate or previously deposited layers, and then an electron beam selectively melts and fuses the powder particles according to the desired part model data. This layer-by-layer process continues until the complete part is built. When the electron beam scans across the powder surface, energy transfers from the top surface to the subsurface. The powder melts when the temperature reaches the material's liquidus temperature. Figure~\ref{fig: EBM} shows the schematic of the electron beam melting machine and a part being built, including the molten pool, solidification zone, substrate, and unmelted powder.

\subsection{Residual stress and treatment} \label{section: bg_residual_stress}


In the electron beam melting process, rapid heating of the upper surface combined with slow heat conduction of the underlying layers creates a steep temperature gradient. This gradient generates convection currents that affect material flow, creates surface tension variations in the molten pool, and induces high internal stresses. Expansion of the heated layer is restricted by the vicinity area, resulting in compressive stress at the surface. Cooling takes place as the heat source moves away, and the contraction of the top layer is then restricted by the surrounding area. This results in tensile residual stress on the top surface. These tensile residual stresses accumulate across multiple scan paths, potentially leading to delamination or crack formation. In addition, the differential cooling and shrinkage between layers causes the top layers to become shorter than the bottom layers. This dimensional mismatch creates a bending moment that distorts the part upward toward the beam direction~\cite{fu20143,li2018residual}.


To mitigate residual stress, in-situ control and post-process control methods have been studied. In-situ controls for additive manufacturing includes both mechanical and thermal methods. Mechanical approaches, such as laser shock peening and rolling, generate compressive pressure to counterbalance residual stresses. However, laser shock peening, while effective, significantly extends printing time and poses integration challenges with powder bed fusion systems~\cite{kalentics2017tailoring}. Thermal gradient control methods reduce residual stress by homogenizing temperature distribution throughout the process. However, they can cause larger grain sizes and anisotropic microstructure, require precise temperature calibration based on alloy properties, and consume more energy~\cite{everton2016review,chia2022process}. Post-heat treatment processes are widely utilized to restore homogeneous and stable microstructures, thereby achieving desired mechanical properties~\cite{vilaro2011fabricated}. However, prior experiments~\cite{shiomi2004residual, mur1996influence} are inconclusive in their optimal treatment temperatures. It is therefore critical to minimize residual stress \textit{during} the fabrication process itself through careful optimization of processing parameters rather than relying solely on post-processing treatments. We discuss the thermal and mechanical models governing the manufacturing process next.

\subsection{Thermal model} \label{section: bg_thermal_model}

The governing equation for the heat transfer analysis in a powder bed fusion process \cite{zinoviev2016evolution} is:
    \begin{equation}
        \rho C_p(T)\frac{\partial T}{\partial t}=\nabla\cdot{(\kappa(T)\nabla T)}+Q_e,
        \label{eq: HeatTrasfer}
    \end{equation}
    where $T(x,y,z,t)$ is the temperature, $\rho$ is the constant density of the material, $C_p(T)$ is the specific heat, $\kappa(T)$ is thermal conductivity, and $Q_e(x,y,z,t)$ is the the applied heat flux. The specific heat and thermal conductivity are modeled with second-degree polynomials as $C_p(T) = a_{0}+a_{1}T+a_{2}T^2$ and $\kappa(T) = b_{0}+b_{1}T+b_{2}T^2$, respectively. The parameters are specified in Section~\ref{section: simulation_runs}. The heat flux due to the electron beam, $Q_e$, can be modeled in a Gaussian form \cite{fu20143, vastola2016controlling}:
    \begin{equation}
        Q_e(x,y,z,t; P, v, r, z_0) = \frac{2AP}{\pi r^2 z_0}\exp{\left( -\frac{2 \left( (x-vt)^2 + y^2 \right)}{r^2} \right)} \frac{1}{5} \left[ -3 \left( \frac{z}{z_0} \right)^2 - 2 \frac{z}{z_0} + 5 \right],
    \end{equation}
    where $P$ is the beam power, $v$ is the scanning speed, $r$ is the beam spot radius, $z_0$ is the beam penetration depth and $A$ is the dimensionless powder absorptivity. The temperature history during the manufacturing process is used as the input to the mechanical model, as well as to obtain $T_{\max}$ in the constraint in \eqref{eq: optimization_setup} to ensure that the metal powder melts.

\subsection{Mechanical model} \label{section: bg_mechanical_model}

The governing equation for mechanical analysis in a powder bed fusion process \cite{megahed2016metal} is given by:
        \begin{equation}
            \nabla\cdot \boldsymbol{\sigma}+\bf=\boldsymbol{0},
        \end{equation}
    where $\boldsymbol{\sigma}(x,y,z,t)$ denotes the stress tensor and $\bf(x,y,z,t)$ are the internal forces that balance external forces. For the electron beam melting process, we first perform the thermal analysis as described above and the resulting temperature fields at every time instance drive the mechanical analysis to calculate stresses and strains. We perform mechanical analysis at every time instance. 

The total strain is decomposed into elastic strain, plastic strain, and thermal strain, $\boldsymbol{\epsilon} = \boldsymbol{\epsilon}^e+\boldsymbol{\epsilon}^p+\boldsymbol{\epsilon}^t$. The elastic stress $\boldsymbol{\sigma}^{e}$ can be related to the elastic strain $\boldsymbol{\epsilon}^{e}$ through the rank-four elasticity tensor~$\boldsymbol{C}$:
        \begin{equation}
            \boldsymbol{\sigma}^{e} = \boldsymbol{C} \boldsymbol{\epsilon}^{e}.
        \end{equation}
The plastic strain $\boldsymbol{\epsilon}^{p}$ is modeled by considering an elastic-perfectly plastic \cite{zhao2015numerical} condition in the model. The thermal strain is calculated from the thermal expansion constitutive relationship: $\boldsymbol{\epsilon}^t = \alpha_t \Delta T$ where $\alpha_t$ is the thermal expansion coefficient. The presence of the thermal strain tensor ensures correct distortion calculation during the melting stage as well as the thermal shrinkage during the cooling stage~\cite{megahed2016metal}.
 
We choose the von Mises stress as the residual stress in the AM part. The von Mises stress is a function of the normal and shear stresses, and is used to predict yielding of materials under complex loading. When the yield stress of the material is reached in the powder bed fusion process, plastic distortion occurs and the quality of the manufactured part is compromised. 

In summary, we use the temperature history from the thermal model to compute the stress and strain in the mechanical model. The coupling between the two models is therefore uni-directional, as the mechanical response depends on the thermal response but not vice versa \cite{debroy2017building}. 

\subsection{Finite element model} \label{section: FEM}

\begin{figure}
    \small
    \captionsetup[subfigure]{font=small}
    \begin{minipage}{0.45\textwidth}
        \centering
        \subcaptionbox{Part geometry and the $x^c-z^c$ plane passing through the part centroid. \label{fig: Part}}{
            \includegraphics[width=\textwidth]{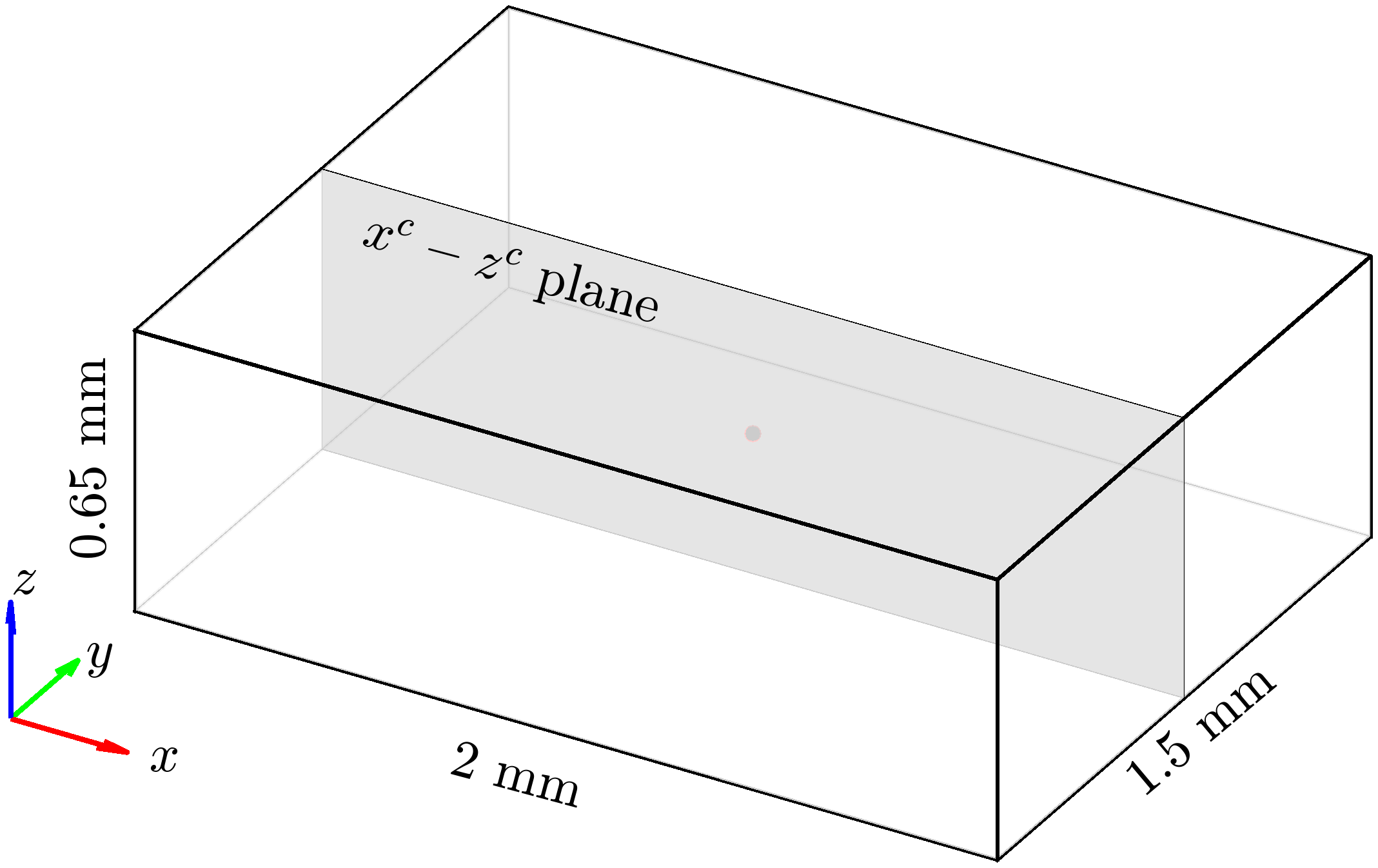}}
    \end{minipage}
    \hfill 
    \begin{minipage}{0.45\textwidth}
        \centering
        \subcaptionbox{Integration points of C3D8R elements on the $x^c-z^c$ plane. We use the maximum von Mises stress in this plane at the end of the cooling step to obtain $\sigma_{\max}$, the residual stress of interest. \label{fig: StressGrid}}{
            \includegraphics[width=\textwidth]{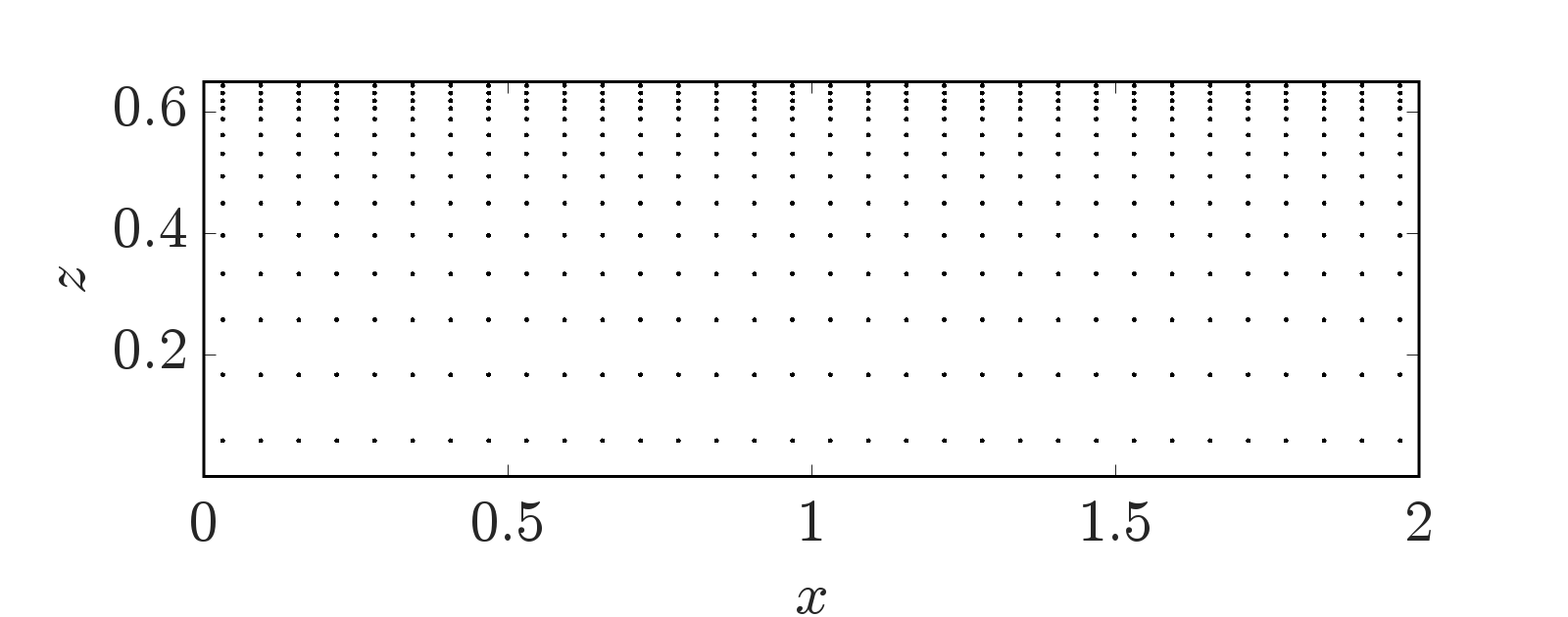}}
    \end{minipage}
    \caption{Part geometry and integration points on the $x^c-z^c$ plane.}
\end{figure}

The AM part of interest is shown in Figure \ref{fig: Part}. The part's length is $l = $ 2~mm, and the width and height are 1.5~mm and 0.65~mm, respectively. We use Ti-6Al-4V powder to manufacture this part. This alloy is commonly used in the aerospace and bioengineering industries because of its good mechanical properties, low density, and good corrosion resistance \cite{chastand2018comparative}.

We simulate the manufacturing of this part via an electron beam melting process in \textsf{Abaqus} using both thermal and mechanical finite element models. In electron beam melting, fused layers are stacked on top of each other, therefore, the layer thickness influences the resolution of the build \cite{cansizoglu2008applications}. Early versions of electron beam melting process equipment used 100 $\mu$m as the standard layer thickness. The current standard layer thickness has been reduced to 50-70~$\mu$m \cite{karlsson2013characterization}. We use layer thickness of 50~$\mu$m here. In order to mitigate the computational cost of the finite element analysis, we only simulate the process of fusing the last powder layer (50~$\mu$m thick) to the bulk material formed by previous scans (0.6~mm high). We consider a single scan along the length in the $x$ direction through the powder layer on the top. 

The inputs to the finite element model include the scanning speed $v$ [mm/s], beam power $P$ [W], preheating temperatures $T_0$ [$^{\circ}$C], yield strength $Y$ [MPa], elastic modulus $E$ [GPa], bulk density $\rho$ [kg/$\text{m}^3$], specific heat $C_p$ [J/(kg$\cdot$K)] and thermal conductivity $\kappa$ [W/(m$\cdot$K)]. Among these inputs, $C_p(T)$ and $\kappa(T)$ are functions of temperature and parameterized by polynomial coefficients, $a_{i}$ and $b_{i}, \ i = 1,2,3$, respectively, as discussed in detail in Section~\ref{section: simulation_runs}. All values are converted to mm, MPa, and related units before being input into \textsf{Abaqus}. 

We use three-node linear brick-type elements, specifically, DC3D8 (one degree of freedom per element, temperature) for the thermal model and C3D8R (24 degrees of freedom per element, $x,y,\ \text{and}\ z$ displacements) for the mechanical model. We use a non-uniform mesh, locally refined for the powder region where the heat flux is applied. We gradually coarsen the mesh for the rest of the part, which is indicated by the integration points corresponding to C3D8R elements intersecting the $x^c-z^c$ plane (the $x-z$ plane passing through the part centroid) in Figure~\ref{fig: Part}. The mesh consists of 13,200 nodes and 10,752 elements in total. At each time instance, the thermal model has 13,200 DoFs and the mechanical model has 39,600 DoFs. 

In the thermal model, we allow the top surface to exchange heat with the surrounding through radiation. We model the radiation using the Stefan-Boltzmann law, $Q_r(x,y,z = 0.65\ \text{mm}, t) = \sigma_{\textsf{SB}} \epsilon_s (T^4(x,y,z=0.65\ \text{mm},t) - T_c^4)$, where $\sigma_{\textsf{SB}} = 5.67 \times 10^{-8}\ \text{W}\cdot \text{m}^{-2}\cdot \text{K}^{-4}$ is the Stefan-Boltzmann constant, ${\epsilon}_s$ is the surface emissivity, and $T_c$ is the constant chamber temperature. In the mechanical model, the boundary surfaces in the $x$-direction and $y$-direction are fixed in the $x$-coordinates and $y$-coordinates respectively. The bottom surface is considered fixed in all coordinates.

There are four steps in the \textsf{Abaqus} models, namely laying powder, preheating, moving heat, and cooling. We simulate the process of laying the new powder on bulk material by activating the initially deactivated elements representing the powder layer. For our optimization problem~\eqref{eq: optimization_setup}, we extract $T_{\max}$ from the temperature history in step 3 (moving heat) and $\sigma_{\max}$ from the residual stress at the end of step 4 (cooling).

\section{Solution to the risk-based design optimization problem} \label{section: DUU_AM}

In this section, we present an efficient solution of the optimization problem formulated in Section~\ref{section: problem_discription}. It is computationally prohibitive to obtain temperature history and residual stress for solving the design problem \eqref{eq: optimization_setup} using the high-fidelity finite element models described in Section~\ref{section: bg_thermal_model} and Section~\ref{section: bg_mechanical_model}, as the optimization process requires repeated model evaluations. To address this computational challenge, we construct surrogate models with training data obtained from the high-fidelity finite element models and solve the design optimization problem with surrogate models.

In Section~\ref{section: risk_measures}, we introduce the risk measure, buffered probability of failure, used as a constraint in the design problem. Then, we discuss the surrogate modeling techniques for both thermal and mechanical finite element models in Section~\ref{section: surrogate_model}. In Section~\ref{section: opt_sol}, we detail the solution of the design optimization problem with the integration of BPOF and surrogate models.
\subsection{Buffered probability of failure (BPOF)} \label{section: risk_measures} 


For an engineering system to be designed under uncertainty, we separate the inputs to the system into design variables $\bd \in \cD \subseteq \real^{N_d}$, random variables $Z \in \Omega$, and parameters $\boldsymbol{\theta}\in \Theta \subseteq \real^{N_{\theta}}$. Here, $\cD$ is the design space, $\Omega$ is the space of random variables, $\Theta$ is the parameter space and $\real$ represents the real numbers. A realization of $Z$ is denoted by $\boldsymbol{z} \in \real^{N_z}$. 

We first introduce concepts related to buffered probability of failure, including probability of failure and superquantile. The failure of an engineering system can be described by the condition $g(\bd, Z; \boldsymbol{\theta}) > \tau$, where $g: \cD \times \Omega \times \Theta \rightarrow \real$ is the limit state function and $\tau \in \real$ is the failure threshold. For a system under uncertainty, $g(\bd, Z;\boldsymbol{\theta})$ is a random variable given a particular design $\bd$ and parameter $\boldsymbol{\theta}$. The probability of failure is defined as 
        \begin{equation}
            p_{\tau}(g(\bd, Z; \boldsymbol{\theta})):= \mathbb{P}[g(\bd, Z; \boldsymbol{\theta}) > \tau].
            \label{eq: POF}
        \end{equation}
    The POF is a measure of the set $\{g(\bd, Z;\boldsymbol{\theta}) > \tau\}$ thus it does not include information about the magnitude of failure.

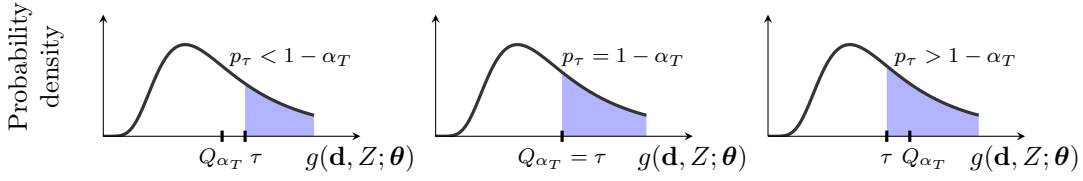
\begin{figure}
    \centering
    \input{fig_Q_alpha_POF_all}
    \caption{Relationship between the probability of failure (shaded area under the pdf curve) and $\alpha$-quantile.}
    \label{fig: Q_alpha_POF}
\end{figure}

The $\alpha$-quantile, $Q_{\alpha}$, also known as the value-at-risk at level $\alpha$, of the random variable $g(\bd, Z;\boldsymbol{\theta})$ can be expressed as:
        \begin{equation}
            Q_{\alpha}[g(\bd, Z;\boldsymbol{\theta})]:=F^{-1}_{g(\bd, Z;\boldsymbol{\theta})}(\alpha),
        \end{equation}
    where $F^{-1}_{g(\bd, Z;\boldsymbol{\theta})}(\cdot)$ is the inverse cumulative density function. Given a target reliability $\alpha_T$, the corresponding target probability of failure is $1 - \alpha_T$, thus, POF and $Q_{\alpha_T}$ are related by:
        \begin{equation}
            \begin{aligned}
                & Q_{\alpha_T}[g(\bd, Z;\boldsymbol{\theta})] < \tau \qquad \Leftrightarrow \qquad p_{\tau}(g(\bd, Z; \boldsymbol{\theta})) < 1-\alpha_T, \\
                & Q_{\alpha_T}[g(\bd, Z;\boldsymbol{\theta})] = \tau \qquad \Leftrightarrow \qquad p_{\tau}(g(\bd, Z; \boldsymbol{\theta})) = 1-\alpha_T, \\
                & Q_{\alpha_T}[g(\bd, Z;\boldsymbol{\theta})] > \tau \qquad \Leftrightarrow \qquad p_{\tau}(g(\bd, Z; \boldsymbol{\theta})) > 1-\alpha_T,
            \end{aligned}
        \end{equation}
    as shown in Figure~\ref{fig: Q_alpha_POF}.

The $\alpha$-superquantile, $\bar{Q}_{\alpha}$, also know as conditional value-at-risk at level $\alpha$, is defined based on the $\alpha$-quantile $Q_{\alpha}$:
        \begin{equation}
            \bar{Q}_{\alpha}[g(\bd, Z;\boldsymbol{\theta})]:= Q_{\alpha}[g(\bd, Z;\boldsymbol{\theta})] + \frac{1}{1-\alpha}\mathbb{E}\Bigr[\bigr[g(\bd, Z;\boldsymbol{\theta}) - Q_{\alpha}[g(\bd, Z;\boldsymbol{\theta})]\bigr]^{+}\Bigr],
            \label{eq: superquantile_definition}
        \end{equation}
    where $[c]^{+} := \max\{0, c\}$. The conditional value-at-risk is a tail expectation, i.e., an average over the portion exceeding the failure threshold. It is the sum of the $\alpha$-quantile and a non-negative term and is thus conservative compared to the $\alpha$-quantile. By definition, $\bar{Q}_0[g(\bd, Z;\boldsymbol{\theta})] = \mathbb{E}[g(\bd, Z;\boldsymbol{\theta})]$ and $\bar{Q}_1[g(\bd, Z;\boldsymbol{\theta})]$ is the essential supremum of $g(\bd, Z;\boldsymbol{\theta})$.


The buffered probability of failure is an alternative measure of reliability that adds a buffer to the POF. The BPOF is defined based on the superquantile as
        \begin{equation}
            \bar{p}_{\tau}(g(\bd,Z;\boldsymbol{\theta})):=
                \begin{cases}
                \{1-\alpha\ |\ \bar{Q}_\alpha[g(\bd,Z;\boldsymbol{\theta})]=\tau\}, & \text{if } \bar{Q}_0[g(\bd,Z;\boldsymbol{\theta})]< \tau <\bar{Q}_1[g(\bd,Z;\boldsymbol{\theta})], \\
                0, & \text{if } \tau \geq \bar{Q}_1[g(\bd,Z;\boldsymbol{\theta})], \\
                1, & \text{otherwise.}
                \end{cases}
            \label{eq: BPOF_definition}
        \end{equation}
Note that by definition, $\mathbb{P}\bigr[g(\bd,Z;\boldsymbol{\theta}) \geq Q_{\alpha}[g(\bd,Z;\boldsymbol{\theta})]\bigr] = 1 - \alpha$, thus the first condition in \eqref{eq: BPOF_definition} is 
        \begin{equation}
            \bar{p}_{\tau}(g(\bd,Z;\boldsymbol{\theta})) = \mathbb{P}\bigr[g(\bd,Z;\boldsymbol{\theta}) \geq Q_{\alpha}[g(\bd,Z;\boldsymbol{\theta})]\bigr],
            \label{eq: BPOF_quantile}
        \end{equation}
    where $\alpha$ satisfies $\bar{Q}_{\alpha}[g(\bd,Z;\boldsymbol{\theta})] = \tau$. The BPOF can be viewed as the probability of exceeding a quantile given the condition on $\alpha$. Alternatively, the BPOF can be written in the form of an expectation as 
    \begin{equation}
        \bar{p}_{\tau}(g(\bd,Z;\boldsymbol{\theta})) = \min_{\zeta< \tau} \frac{\mathbb{E}\bigr[[g(\bd,Z;\boldsymbol{\theta}) - \zeta]^+\bigr]}{\tau - \zeta},
    \end{equation}
    where $\zeta$ is an auxiliary variable, and the optimum $\zeta^*$ is the threshold from \eqref{eq: BPOF_quantile} that provides $\mathbb{E}[g(\bd,Z;\boldsymbol{\theta})|g(\bd,Z;\boldsymbol{\theta}) > \zeta^*] = \tau$. This form of BPOF assumes that $\bar{Q}_0[g(\bd, Z;\boldsymbol{\theta})] < \tau < \bar{Q}_1[g(\bd, Z;\boldsymbol{\theta})]$ and $g(\bd, Z;\boldsymbol{\theta})$ is integrable. 
    
\begin{figure}[!htbp]
    \centering
    \input{fig_BPOF_buffer}
    \caption{BPOF $\bar{p}_{\tau}$ illustration at threshold $\tau$: BPOF equals POF plus the buffer.}
    \label{fig: BPOF_buffer}
\end{figure}
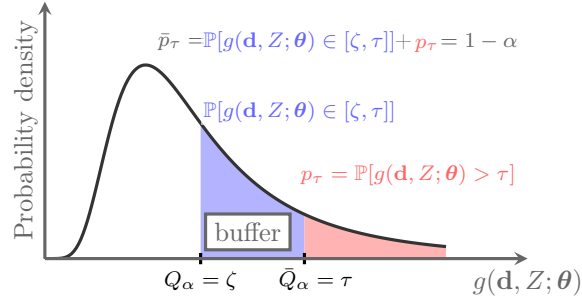

The conservativeness of BPOF stems from the selection of the threshold $\zeta \leq \tau$ based on data. When realizations of $g(\bd, Z;\boldsymbol{\theta})$ beyond $\tau$ are large (corresponding to potentially catastrophic failures), $\zeta$ has to be smaller to drive the expectation beyond $\zeta$ to $\tau$ thus making BPOF bigger. The $\zeta$-tail also has to have an average equal to $\tau$, therefore, when there are a large number of near-failure events, the near-failure region $[\zeta, \tau]$ is determined by the frequency and magnitude of tail events around $\tau$. This can be intuitively viewed as a buffer to the POF as shown in Figure~\ref{fig: BPOF_buffer}. 

The BPOF in \eqref{eq: BPOF_quantile} can be further split into a sum of a non-negative term and POF as 
        \begin{equation}
            \begin{aligned}
                \bar{p}_{\tau}(g(\bd,Z;\boldsymbol{\theta})) &= \mathbb{P}\bigr[g(\bd,Z;\boldsymbol{\theta}) \in [Q_{\alpha}[g(\bd,Z;\boldsymbol{\theta})], \tau]\bigr] + \mathbb{P}[g(\bd,Z;\boldsymbol{\theta}) > \tau]\\
                &= \mathbb{P}\bigr[g(\bd,Z;\boldsymbol{\theta}) \in [\zeta, \tau]\bigr] + p_{\tau}(g(\bd,Z;\boldsymbol{\theta})),
                \label{eq: BPOF_POF}
            \end{aligned}
        \end{equation}
    where $\zeta = Q_{\alpha}[g(\bd,Z;\boldsymbol{\theta})]$ is affected by the condition on $\alpha$ through the superquantile. In summary, BPOF considers the magnitude of the failure in addition to failure frequency through the first term in \eqref{eq: BPOF_POF} and is a conservative estimate of the POF. 

We can use a sampling-based approach, see Algorithm~\ref{alg: MC_estimation}, to estimate the $\alpha$-quantile $Q_{\alpha}$, $\alpha$-superquantile $\bar{Q}_{\alpha}$, and BPOF $\bar{p}_{\tau}$ as $Q_{\alpha}$ and $\bar{Q}_{\alpha}$ can be viewed as expectations. With $m$ Monte Carlo samples, the estimation errors decrease at a rate of $1/\sqrt{m}$ and the estimator variance increases with larger $\alpha$ values.

\begin{algorithm}
\caption{Sampling-based estimations of $Q_{\alpha},\ \bar{Q}_{\alpha},$ and $\bar{p}_{\tau}$}
\hspace*{\algorithmicindent} \textbf{Input} $m$ i.i.d. samples $\boldsymbol{z}_1, \ldots, \boldsymbol{z}_m$ of random variable $Z$; design variable $\bd$; parameters $\boldsymbol{\theta}$; risk level $\alpha \in (0, 1)$.\\
\hspace*{\algorithmicindent} \textbf{Output} Sample approximations $\widehat{Q}_\alpha[g(\bd,Z;\boldsymbol{\theta})]$, $\widehat{\bar{Q}}_\alpha[g(\bd,Z;\boldsymbol{\theta})]$, and $\widehat{p}_{\tau}(g(\bd,Z;\boldsymbol{\theta}))$.
    \begin{algorithmic}[1]
    \State Evaluate limit state function at the samples to get $g(\bd,\boldsymbol{z}_1;\boldsymbol{\theta}), \ldots, g(\bd,\boldsymbol{z}_m;\boldsymbol{\theta})$.
    \State Sort values of limit state function in descending order and relabel the samples so that
    $$g(\bd,\boldsymbol{z}_1;\boldsymbol{\theta}) > g(\bd,\boldsymbol{z}_2;\boldsymbol{\theta}) > \cdots > g(\bd,\boldsymbol{z}_m;\boldsymbol{\theta}).$$
    
    %
    \State $\widehat{Q}_{\alpha}[g(\bd,Z;\boldsymbol{\theta})] \leftarrow g(\bd, \boldsymbol{z}_{k_\alpha};\boldsymbol{\theta})$, where $k_\alpha = m(1 - \alpha)$.
    \State Estimate $\widehat{\bar{Q}}_{\alpha}[g(\bd,Z;\boldsymbol{\theta})]$ using \eqref{eq: superquantile_definition}.
    %
    %
    \State $c = g(\bd, \boldsymbol{z}_1;\boldsymbol{\theta})$. 
    \State $k = 1$.
    \State \textbf{while} $c \geq \widehat{\bar{Q}}_{\alpha}[g(\bd,Z;\boldsymbol{\theta})]$ \textbf{do} 
    \State \quad $k \leftarrow k + 1$.
    \State \quad $c = \frac{1}{k}\sum_{i=1}^k g(\bd, \boldsymbol{z}_k; \boldsymbol{\theta})$. 
    \State \textbf{end while}
    \State Estimate BPOF as $\widehat{p}_{\tau}(g(\bd, Z; \boldsymbol{\theta})) \approx \frac{k-1}{m}$, $\tau \approx c$. 
\end{algorithmic}\label{alg: MC_estimation}
\end{algorithm}
\subsection{Surrogate modeling techniques} \label{section: surrogate_model}
We use surrogate models to approximate the FEM solutions and accelerate the optimization. Here, we discuss the surrogate modeling techniques for both the thermal and mechanical models.

The temperature history in additive manufacturing processes exhibits strong temporal correlation, while residual stress exhibits spatial correlation in the printed part. We use an efficient approach to transform these correlated high-dimensional outputs into an uncorrelated feature space, which enables separate surrogate modeling of features as discussed below. Specifically, we first map temperature histories and residual stresses to features (see Section~\ref{section: SVD} for details) in uncorrelated spaces using singular value decomposition (SVD). We then identify the active subspaces of the input variables and construct surrogate models of features versus their corresponding active variables in the active subspace. This strategy achieves effective dimension reduction for both the inputs and the QoIs, which enables efficient modeling for the complex physical phenomena inherent in the additive manufacturing process. Finally, we use surrogate model predictions of features and SVD components to recover the temperature histories and residual stress, which are used in the optimization formulation.

\subsubsection{Singular value decomposition for feature generation} \label{section: SVD}

We model the manufacturing process by the design vector $\bd$, the random vector $Z$, and the deterministic parameters $\boldsymbol{\theta}$ as discussed in Section~\ref{section: problem_discription}. To obtain training data for the surrogate models, we  generate $M$ random samples of \textit{both} the design variables and random variables. We denote $\Xi = [v, P, T_0, Y, E, \rho]^{\top}$ as a random variable and a realization of $\Xi$ as $\boldsymbol{\xi}$. We obtain both the thermal and mechanical FEM outputs corresponding to each sample $\boldsymbol{\xi}_i$, $i=1, 2, \ldots, M$, and reshape the FEM outputs corresponding to each sample as a row vector of length $N$. We then assemble these vectors into a matrix $\bF \in \real^{M\times N}$. 

The singular value decomposition of $\bF \in \mathbb{R}^{M \times N}$ is $\bF=\bU \boldsymbol{\Sigma} \bV^{\prime}=\sum_{k=1}^{\min (M, N)} \sigma_{k} \bu_{k} \bv_{k}^{\prime}$, where  $\bU \in \mathbb{R}^{M \times M}$ and $\bV \in \mathbb{R}^{N \times N}$ are unitary matrices that have the left singular vectors $\bu_{k}$ and right singular vector $\bv_{k}$ as columns, respectively. The diagonal matrix $\boldsymbol{\Sigma} \in \mathbb{R}^{M \times N}$ contains singular values $\sigma_1, \sigma_2, \ldots, \sigma_p$, $p = \min\{M,\ N\}$ which are ordered by magnitude, $\sigma_1 \geq \sigma_2 \geq \ldots \geq \sigma_p \geq 0$. 

The best approximation in the least-squares sense of the data matrix $\bF$ can be obtained by $\hat{\bF} = \bU_k\boldsymbol{\Sigma}_k\bV_k'$, where $\bU_k$ is a matrix containing the first $k$ left singular vector, $\boldsymbol{\Sigma}_k$ is the first $k$ singular values organized in a $k\times k$ diagonal matrix, and $\bV_k$ is a matrix containing the first $k$ right singular vectors. We denote the product of the matrices $\bU_k$ and~$\boldsymbol{\Sigma}_k$ as:
        \begin{equation}
            \cF = \bU_k\boldsymbol{\Sigma}_k \in \real^{M\times k},
        \end{equation}
which induces a factorization of $\hat{\bF}$. Due to the orthonormality of $\bV$ in the SVD, $\cF$ has the following properties: (a) each row $\cF_{i, :}$ can be viewed as a $k$-dimensional coordinate on the orthonormal basis $[\bv_1, \bv_2, \ldots, \bv_k]=: \bV[:, 1:k]$. Any two components in the coordinate, $\cF_{i, j}$, and $\cF_{i, m}$, $j \neq m$, are uncorrelated. We refer to $\cF_{i, j}$, the coordinate component in $\cF_{i, :}$ that is indexed by $j$, as a realization of \textit{feature}~$j$. (b) The mapping between $\cF$ and $\hat{\bF}$ is exact with the knowledge of $\bV_k$: 
    \begin{equation}
        \hat{\bF} = \cF\bV_k'.
        \label{eq: feature_to_F}
    \end{equation}
    When $k = N$, this becomes $\bF = \cF\bV'$.  

In summary, the original data in $\bF$ can be appoximated with $k$ features, $\cF_{:, 1}, \cF_{:, 2}, \ldots, \cF_{:, k}$, in $\cF$ and the first $k$ right singular vectors of $\bF$. These features, corresponding to the first $k$ left singular vectors and singular values, are uncorrelated. When $k$ is smaller than $N$, dimension reduction is achieved. 

\subsubsection{Active subspace discovery and active variable calculation} \label{section: AS}

An active subspace is a low-dimensional subspace that consists of important directions in a model's input space~\cite{constantine2015active}. Most of the variability in the output due to the uncertain inputs is captured along these important directions. In the design optimization problem, we are interested in discovering the active subspace related to the features corresponding to the thermal and mechanical model outputs. Every feature column $\cF_{:, j}$ can be treated as a scalar-valued function of the input $\boldsymbol{\xi}$, $\cF_{:,j}(\boldsymbol{\xi})$. An active subspace is thus a low-dimensional subspace in the input domain that effectively captures the variability in $\cF_{:, j}$ due to variations in $\boldsymbol{\xi}$. 

The directions constituting the active subspace are the dominant eigenvectors of the covariance matrix $\cC=\int_{\Omega}\left(\nabla_{\Xi} \cF_{:, j} \right)\left(\nabla_{\Xi} \cF_{:, j} \right)^{\top} \pi(\boldsymbol{\xi})\text{d} \boldsymbol{\xi}$, where $\pi(\boldsymbol{\xi})$ is the joint probability density function of ${\Xi}$. The eigenvalue decomposition of the matrix $\cC$ is given by $\cC=\bW\boldsymbol{\Lambda}\bW^{'}$. Since $\cC$ is real, symmetric and positive semidefinite, its eigenvalues are non-negative real values and its eigenvectors are orthogonal. Thus, the eigenvector matrix and eigenvalue matrix can be partitioned about the $r$th eigenvalue such that there is a significant drop in magnitude, i.e., $\lambda_r/\lambda_{r+1} \gg 1$ as 
    \begin{equation}
        \bW=\left[\bW_{1}\  \bW_{2}\right], \quad \mathbf{\Lambda}=\left[\begin{array}{cc}
        \mathbf{\Lambda}_{1} & \\
        & \mathbf{\Lambda}_{2}
        \end{array}\right].
        \label{eq: AS_W_Lambda}
    \end{equation}
    The columns of $\bW_{1} = [\bw_1, \bw_2, \ldots, \bw_r]$ span the dominant eigenspace of $\cC$ and the active subspace and the corresponding active variable is calculated as $\boldsymbol{\eta} = \bW_1^{'} \boldsymbol{\xi} \in \real^r$. When $r$ is smaller than the dimension of $\Xi$, dimension reduction is achieved.

We use the trained surrogate models to estimate feature values for unseen model inputs, including design variables, random variables, and system parameters. We then use the right singular vector matrix $\bV$ of the SVD to obtain the temperature history and residual stress estimations corresponding to the unseen input using \eqref{eq: feature_to_F}. These estimations are used in the optimization problem~\eqref{eq: optimization_setup}. We only compute surrogate models for $K_{\cT}$ features of the temperature history and $K_{\cS}$ features of the residual stress. We discuss the details of surrogate model constructions in Section~\ref{section: soln_surr_mdl_ther} and Section~\ref{section: soln_surr_mdl_mech}.

\subsubsection{Surrogate model for temperature history} \label{section: soln_surr_mdl_ther}

In the optimization setup \eqref{eq: optimization_setup}, the constraint on temperature history ensures the powder properly melts. Since we model a single scan along the length at the top surface of the part as discussed in Section~\ref{section: background}, we choose to record the temperature history at the centroid of the top surface of the part. The total scan time $t_{\text{scan}} = l/v$ varies between simulations with different scanning speed.

Since the beam moves at constant velocity, the temperature is expected to peak at approximately $0.5 \cdot T_{\text{scan}}$. We take temperature measurements at selected non-uniform time instances within $[0, t_{\text{scan}}]$ according to the following distribution: (a) ten uniformly spaced measurements, $T_1,T_2,\ldots,T_{10}$, in the interval $[0, 0.405 \cdot T_{\text{scan}}]$; (b) ten uniformly spaced measurements, $T_{11},T_{12},\ldots,T_{20}$, in the interval $[0.45 \cdot T_{\text{scan}}, 0.54 \cdot T_{\text{scan}}]$; (c) eleven uniformly spaced measurements, $T_{21},T_{22},\ldots,T_{31}$, in the interval $[0.55 \cdot T_{\text{scan}}, T_{\text{scan}}]$. Note that we use an increased sampling density within $[0.45 \cdot T_{\text{scan}}, 0.54 \cdot T_{\text{scan}}]$ to ensure precise capture of the maximum temperature at the center of the part's top surface.

With $\cM$ simulations, we compile a snapshot matrix $\bT \in \real^{\cM\times 31}$. Each row in the snapshot matrix, $\bT_{i, :}$, contains the temperatures at thirty-one time instances $[T_1,T_2,\ldots,T_{31}]$ of simulation $i$, where $i=1,2,\ldots,\cM$. We apply an SVD and active subspace discovery, as detailed in Sections~\ref{section: SVD}~and~\ref{section: AS}, to extract features $\cT$ and corresponding active variables $\boldsymbol{\eta}_{(\bT)}$. We select the most significant $K_{\cT}$ features to construct surrogate models $\hat{G}_{\cT_{:,j}}(\boldsymbol{\eta}_{(j_{\bT})}) \approx \cT_{:,j}((\bd, Z;\boldsymbol{\theta}))$ for $j = 1, 2, \ldots, K_{\cT}$. We choose $K_{\cT}$ based on prediction error on the training data, see Section~\ref{section: surr_mdl_ther} for a detailed discussion.

\subsubsection{Surrogate model for residual stress} \label{section: soln_surr_mdl_mech}

We measure the residual stress on the non-uniform grid comprising 32 points along the length and 14 points along the height on the $x^c-z^c$ plane, as shown in Figure~\ref{fig: StressGrid}. We use the maximum von Mises stress measured at the end of the cooling step in this plane as the residual stress of interest. 

Each simulation generates a comprehensive stress field $\boldsymbol{\sigma} \in \mathbb{R}^{32\times14}$ for the $x^c-z^c$ plane. For a set of $\cM$ simulations, we vectorize the stress fields and consolidate them into a matrix $\bS \in \mathbb{R}^{\cM\times 448}$. We then apply an SVD and active subspace discovery, as discussed in Sections~\ref{section: SVD}~and~\ref{section: AS}, to obtain features $\cS$ and corresponding active variables $\boldsymbol{\eta}_{(\bS)}$. We select the most influential $K_{\cS}$ features to construct surrogate models $\hat{G}_{\cS_{:,j}}(\boldsymbol{\eta}_{(j_{\bS})}) \approx \cS_{:,j}((\bd, Z;\boldsymbol{\theta}))$ for $j = 1, 2, ..., K_{\cS}$. Similar to the thermal case discussed above, we choose $K_{\cS}$ based on prediction error on the training data, see Section~\ref{section: surr_mdl_mech} for a detailed discussion.

\subsection{Optimization problem solution} \label{section: opt_sol}

The original formulation in \eqref{eq: optimization_setup} can be reformulated to an optimization problem involving an expectation in the constraint on residual stress:
    \begin{equation}
        \begin{aligned}
        & \underset{\substack{P_L \leq P \leq P_U, \\v_L \leq v \leq v_U,\\ \zeta < \tau}}{\text{min}}
        & & P\cdot\frac{l}{v} \\
        & \text{subject to}
        & & \frac{\mathbb{E}[\sigma_{\max}(\bd, Z; \boldsymbol{\theta}) - \zeta]^{+}}{\tau-\zeta} \leq 1 - \alpha_T,\\
        & & & T_{\text{liq}} < T_{\max}(\bd, Z; \boldsymbol{\theta}) < 1.1 \times T_{\text{liq}},
        \end{aligned}
        \label{eq: optimization_reformulated}
    \end{equation}
    where $\zeta$ is an auxiliary variable and $\tau$ is the failure threshold related to the reliability level $\alpha$ as discussed in Section~\ref{section: risk_measures}. The optimization formulation in \eqref{eq: optimization_reformulated} can be efficiently implemented using the surrogate models of the temperature history and residual stress discussed above, as detailed in Algorithm~\ref{alg: solution}.

\begin{algorithm}
\caption{Solution to the risk-based design optimization~\eqref{eq: optimization_reformulated}}
    \begin{algorithmic}[1]
    \Algphase{Phase 1 - Surrogate model constructions}
    \State Generate $\cM$ samples of $\bd$, $\boldsymbol{d}_1,\ldots, \boldsymbol{d}_{\cM}$, and realizations of $Z$, $\boldsymbol{z}_1,\ldots, \boldsymbol{z}_{\cM}$; denote $\boldsymbol{\xi}_i = [\boldsymbol{d}_i^{\top}, \boldsymbol{z}_i^{\top}]^{\top}$.
    \State Perform FEM simulations using $[\boldsymbol{\xi}_i^{\top}; \boldsymbol{\theta}^{\top}]^{\top}, i = 1,2,\ldots, \cM$, as inputs and obtain the temperature snapshot matrix $\bT$ and residual stress matrix $\bS$.
    \State Perform an SVD $\bT = \bU_{\bT} \boldsymbol{\Sigma}_{\bT} \bV_{\bT}^{\prime}$.
    \State Select $K_{\cT}$ features, $\cT_{:, 1}, \cT_{:, 2},\ldots,\cT_{:, K_{\cT}}$.
    \State \textbf{for} $j_{\bT} = 1, 2, \ldots, K_{\cT}$  \textbf{do}
    \State \quad Perform an active subspace discovery to find the matrix $\bW_{1(j_{\bT})}$.
    \State \quad Construct surrogate model $\hat{G}_{\cT_{:, j_{\bT}}}$ with training data set $\left\{\left(\bW_{1(j_{\bT})}^{\top}\cdot\boldsymbol{\xi}_i,\ \cT_{i, j_{\bT}}\right)\right\}_{i = 1}^{\cM}$.
    \State \textbf{end for}
    \State Perform an SVD $\bS= \bU_{\bS} \boldsymbol{\Sigma}_{\bS} \bV_{\bS}^{\prime}$.
    \State Select $K_{\cS}$ features, $\cS_{:, 1}, \cS_{:, 2},\ldots,\cS_{:, K_{\cS}}$.
    \State \textbf{for} $j_{\bS} = 1, 2, \ldots, K_{\cS}$  \textbf{do}
    \State \quad Perform an active subspace discovery to find the matrix $\bW_{1(j_{\bS})}$.
    \State \quad Construct surrogate model $\hat{G}_{\cS_{:, j_{\bT}}}$ with training data set $\left\{\left(\bW_{1(j_{\bS})}^{\top}\cdot\boldsymbol{\xi}_i,\ \cS_{i, j_{\bS}}\right)\right\}_{i = 1}^{\cM}$.
    \State \textbf{end for}
    \Algphase{Phase 2 - Optimization iterations}
    \State \textbf{Initialize}: Set $\bd = \bd^{\{0\}}$, $i=1$.
    \State \textbf{while} not converged \textbf{do}
    \State \quad Generate $\cN$ samples of $Z$, $\boldsymbol{z}_1, \boldsymbol{z}_2, \ldots, \boldsymbol{z}_{\cM}$.
    \State \quad \textbf{for} $k = 1,2,\ldots, \cN$ \textbf{do}
    \State \quad \quad \textbf{for} $j_{\bT} = 1, 2, \ldots, K_{\cT}$  \textbf{do}
    \State \quad \quad \quad Calculate active variable $\boldsymbol{\eta}_{k(j_{\bT})} = \bW_{1(j_{\bT})}^{\top} \cdot \boldsymbol{\xi}_k$, where $\boldsymbol{\xi}_k = [\bd^{\top}, \boldsymbol{z}_j^{\top}]^{\top}$.
    \State \quad \quad \quad Calculate feature value $\hat{G}_{\cT_{:, 1}}(\boldsymbol{\eta}_{k(j_{\bT})})$.
    \State \quad \quad \textbf{end for}
    \State \quad \quad Estimate temperature snapshot $\hat{\bT}_{k, :} \in \real^{1\times31}$, (the $k$-th row of matrix $\hat{\bT}$) as 	
    \begin{equation*}
    \hat{\bT}_{k, :} = [\hat{G}_{\cT_{:, 1}}, \hat{G}_{\cT_{:, 1}},\ldots,\hat{G}_{\cT_{:, K_{\cT}}}]\cdot\bV_{\bT}\left[:, 1:K_{\cT}\right]^{\top}.
    \end{equation*}
    \State \quad \quad Find the maximum temperature $\hat{T}_{k, \max} = \max\{\hat{\bT}_{j, :}\}$.
    \State \quad \quad \textbf{for} $j_{\bS} = 1, 2, \ldots, K_{\cS}$  \textbf{do}
    \State \quad \quad \quad Calculate active variable $\boldsymbol{\eta}_{k(j_{\bS})} = \bW_{1(j_{\bS})}^{\top} \cdot \boldsymbol{\xi}_k$.
    \State \quad \quad \quad Calculate feature value $\hat{G}_{\cS_{:, 1}}(\boldsymbol{\eta}_{k(j_{\bS})})$.
    \State \quad \quad \textbf{end for}
    \State \quad \quad Estimate residual stress $\hat{\bS}_{k, :} \in \real^{1\times448}$, (the $j$-th row of matrix $\hat{\bS}$) as 
    \begin{equation}
    \hat{\bS}_{k, :} = [\hat{G}_{\cS_{:, 1}}, \hat{G}_{\cS_{:, 1}},\ldots,\hat{G}_{\cS_{:, K_{\cS}}}]\cdot\bV_{\bS}\left[:, 1:K_{\cS}\right]^{\top}.
    \label{eq: stress_estimates}
    \end{equation}
    \State \quad \quad Find the maximum residual stress $\hat{\sigma}_{k, \max} = \max\{\hat{\bS}_{k, :}\}$.
    \State \quad Estimate $T_{\max}$ as $\displaystyle \hat{T}_{\max} = \frac{1}{\cN}\sum_{k = 1}^{\cN}\hat{T}_{k, \max}.$
    \State \quad \textbf{end for}
    \State \quad Use the estimated $\hat{T}_{\max}$ and $k$ realizations of $\hat{\sigma}_{\max}$ in~\eqref{eq: optimization_reformulated} in an optimization solver for $\bd^{\{i\}}$.
    \State \quad $\bd \leftarrow \bd^{\{i\}}$, $i \leftarrow i+1$.
    \State \textbf{end while}
    \State \textbf{Return}: Set solution $\bd^* = \bd^{\{i\}}$.
    \end{algorithmic}\label{alg: solution}
\end{algorithm}

\section{Numerical results} \label{section: numerical_results}
%
We present the numerical results in this section. In Section~\ref{section: simulation_runs}, we outline the numerical setup for implementing the finite element simulations. We show surrogate model details and optimization results in Section~\ref{section: nm_surr_mdl} and Section~\ref{section: nm_opt_results}, respectively. Finally, in Section~\ref{section: nm_val}, we demonstrate the validation of the optimized results. 
\subsection{Simulation setup} \label{section: simulation_runs}

The set of process parameters includes scanning speed ($v$), beam power ($P$), and preheating temperature ($T_0$). Mechanical properties include yield strength ($Y$), elastic modulus ($E$), and bulk density ($\rho$). Thermal properties include specific heat ($C_p$) and bulk thermal conductivity ($\kappa$). As discussed in Section~\ref{section: problem_discription}, $\bd = [v, P]^{\top}$ contains the design variables and $Z = [T_0, Y, E, \rho]^{\top}$ is the vector of random variables. The upper and lower bounds of the design variables are $v_L = 100$~mm/s, $v_U = 1{,}000$~mm/s, $P_L = 20$~W, and $P_U = 200$~W. The bounds on the random variables are listed in Table~\ref{tab: rv}. We use uniform distributions for all random variables. The parameter vector $\boldsymbol{\theta}\in \real^6$ consists of six coefficients, $a_{i}$ and $b_{i}, \ i = 1,2,3,$ for calculating the specific heat $C_p(T) = a_{0}+a_{1}T+a_{2}T^2$ and the thermal conductivity $\kappa(T) = b_{0}+b_{1}T+b_{2}T^2$, respectively. We determine the six coefficients using data from \cite{fu20143} and their values are: $a_{0} = 540$~J/(kg$\cdot$K), $a_{1} = 0.43$~J/(kg$\cdot$K$^{2}$), $a_{2} = -3.2\times 10^{-5}$~J/(kg$\cdot$K$^{3}$), $b_{0} = 7.2$~W/(m$\cdot$K), $b_{1} = 0.011$~W/(m$\cdot$K$^{2}$), and $b_{2} = 1.4\times 10^{-6}$~W/(m$\cdot$K$^{3}$). The constant chamber temperature is $T_c$ = 650~${ }^{\circ}$C in the thermal model.

\begin{table}[ht!]
\caption{Uniform random variables used in the additive manufacturing example.} \label{tab: rv}
\centering
{
\begin{tabular}{cccc}
\hline
Random variable  & Lower bound & Upper bound \\
\hline
$T_0$ [${ }^{\circ}$C] & 585 & 715 \\
\hline
$Y$ [MPa] & 742.5 & 907.5 \\
\hline
$E$ [GPa] & 100 & 120 \\
\hline
$\rho$ [kg/m$^3$] & 550.8 & 673.2\\
\hline
\end{tabular}
}
\end{table}

As discussed in Section~\ref{section: FEM}, the finite element simulations use Ti-6Al-4V powder and it is characterized by the following properties \cite{welsch1993materials}. The bulk density is $4{,}428 \mathrm{kg} / \mathrm{m}^3$, powder density is $2{,}700 \mathrm{kg} / \mathrm{m}^3$, solidus temperature is $1{,}605^{\circ} \mathrm{C}$, liquidus temperature is $T_{\text{liq}} = 1{,}650^{\circ} \mathrm{C}$, elastic modulus is $E = 110$~GPa, Poisson's ratio is 0.41, yield strength is $Y = 825$~MPa, surface emissivity is $\epsilon_s = 0.35$, powder absorptivity is $A$ = 0.203, and the thermal expansion coefficient is $\alpha_t = 1\times 10^{-5}$. The latent heat of fusion is $365{,}000 \mathrm{J} / \mathrm{kg}$, which represents the energy required to transition from solid to liquid phase without a change in temperature. To simplify calculation, we use a constant density value $\rho = 4{,}300 \ \mathrm{kg} / \mathrm{m}^3$, the volume-averaged density of powder and bulk alloy. We use a total of 120 training samples ($\cM = 120$ in Algorithm~\ref{alg: solution}) of $(\bd, Z)$, generated using optimal symmetric Latin hypercube sampling in the space $\cD \times \Omega$, to train the surrogate models $\hat{G}_{\cT_{:,j}}(\boldsymbol{\eta}_{(j_{\bT})})$ and $\hat{G}_{\cS_{:,j}}(\boldsymbol{\eta}_{(j_{\bS})})$.

\subsection{Surrogate models} \label{section: nm_surr_mdl}
 
\subsubsection{Temperature snapshot surrogate model} \label{section: surr_mdl_ther}

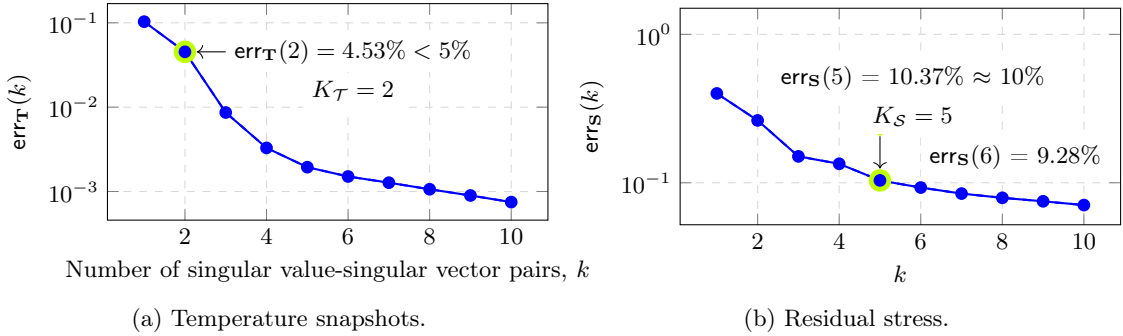
\begin{figure}
    \small
        \captionsetup[subfigure]{font=small}
        \subcaptionbox{Temperature snapshots.\label{fig: T_L_2}}[.48\textwidth]{
        \input{fig_T_L_2}} 
        \captionsetup[subfigure]{font=small}
        \subcaptionbox{Residual stress.\label{fig: S_L_2}}[.48\textwidth]{
        \input{fig_S_L_2}}
    \caption{The mean relative error between (a) temperature snapshots or (b) residual stress and their corresponding approximations using $k$ singular value-singular vector pairs. We choose (a) $K_{\cS}$ or (b) $K_{\cT}$ features for which we build surrogate models based on the error.}
    \label{fig: L_2}
\end{figure}

We perform an SVD on $\bT \in \real^{120\times31}$ and obtain a total of 31 features, $\cT_j, \ j = 1,2, ..., 31$, as discussed in Section~\ref{section: soln_surr_mdl_ther}. We use an error metric $\textsf{err}_{\bT}$ to decide the optimal number of features, $K_{\cT}$. The error metric is defined as:
        \begin{equation}
            \textsf{err}_{\bT}(k) = \frac{1}{\cM}\sum_{i = 1}^{\cM}\frac{\lVert \bT_{i, :} - \hat{\bT}_{i, :}(k)\rVert_2}{\lVert \bT_{i, :}\rVert_2}.
        \end{equation}
It is an averaged relative error between the actual temperature snapshot matrix, $\bT_{i, :}$, and its approximation, $\hat{\bT}_{i, :}(k)$. Note that we estimate $\hat{\bT}_{i, :}(k)$ using the $k$ singular value-singular vector pairs. We calculate the average error with $\cM$ = 120 training samples. We plot $\textsf{err}_{\bT}(k)$ versus $k$ in Figure~\ref{fig: T_L_2}. It can be seen that the top two features yield an error of less than 5\%, therefore, we choose $K_{\cT} = 2$.

\begin{figure}[h!]
    \small
        \captionsetup[subfigure]{font=small}
        \subcaptionbox{Feature 1.}[.48\textwidth]{
        \input{fig_T_surrogate_1}}
        \subcaptionbox{Feature 2.}[.48\textwidth]{
        \input{fig_T_surrogate_2}}
        \caption{Surrogate models for the two features of the temperature data.} 
        \label{fig: T_surrogate}
\end{figure}
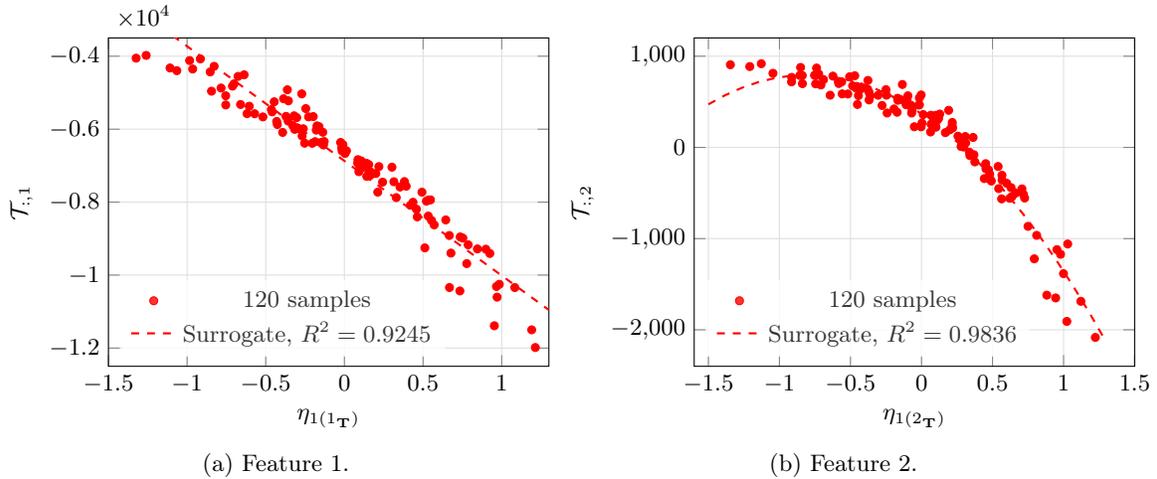

The surrogate models for the top two features are trained with $\cM$=120 samples as described in Algorithm~\ref{alg: solution}. The accuracy of the surrogate model is evaluated by using the coefficient of determination, $R^2 = 1 - \textsf{SS}_{\text{res}}/\textsf{SS}_{\text{tot}}$, where $\textsf{SS}_{\text{res}}$ is the residual sum of squares and $\textsf{SS}_{\text{tot}}$ is the total sum of squares. The $R^2$ score measures the proportion of variance in the model output, and an $R^2$ value close to one indicates a high accuracy. We use a linear regression model as the form as surrogate model for feature 1 and a parabolic function for feature 2. Both features incorporate only one active variable, achieving $R^2$ values exceeding 0.9, as seen in Figure~\ref{fig: T_surrogate}.

\subsubsection{Residual stress surrogate model} \label{section: surr_mdl_mech}


We perform an SVD on $\bS \in \real ^{120\times 448}$  and obtain a total of 120 features, $\cT_j, \ j = 1,2, ..., 120$, as discussed in Section~\ref{section: soln_surr_mdl_mech}. Similar to the thermal case, we use an error metric \textsf{err}$_{\bS}$ to select the optimal number of features, $K_{\cS}$. The error metric is defined as:
        \begin{equation}
            \textsf{err}_{\bS}(k) = \frac{1}{\cM}\sum_{i = 1}^{\cM}\frac{\lVert \bS_{i, :} - \hat{\bS}_{i, :}(k)\rVert_2}{\lVert \bS_{i, :}\rVert_2}.
        \end{equation}
It is an averaged relative error between the actual temperature snapshot matrix, $\bS_{i, :}$, and its approximation, $\hat{\bS}_{i, :}(k)$. Similar to the thermal case, we estimate $\hat{\bS}_{i, :}$ using the $k$ singular value-singular vector pairs and calculate the average error with $\cM$ = 120 training samples.

We plot $\textsf{err}_{\bS}$ versus $k$ in Figure~\ref{fig: S_L_2}. We observe that $\textsf{err}_{\bS}$ decays slower compared to $\textsf{err}_{\bT}$, which can be attributed to the stress calculation's inherent complexity. The stress is derived from the mechanical finite element model using the thermal model's output, as discussed in Section~\ref{section: FEM}. The one-way thermal-mechanical coupling introduces additional computational intricacies. To achieve $\textsf{err}_{\bS} < 5\%$, we need $K_{\cS} > 5$, i.e., the number of features for surrogate modeling becomes prohibitively large. The later features also prove challenging to model accurately, for example, for feature 6, the coefficient of determination $R^2 = 0.3433$ with a sixth-degree polynomial surface using two active variables. In addition, the value of $\textsf{err}_{\bS}$ only decreases from $10.37\%$ to $9.28\%$ when $k$ increases from 5 to 6. Therefore, we choose $K_{\cS} = 5$. 

\begin{figure}[h!]
    \small
        \captionsetup[subfigure]{font=small}
        \subcaptionbox{Feature 1.}[.5\textwidth]{
        \input{fig_S_surrogate_1}}
        \subcaptionbox{Feature 2.}[.5\textwidth]{
        \includegraphics[width=\linewidth]{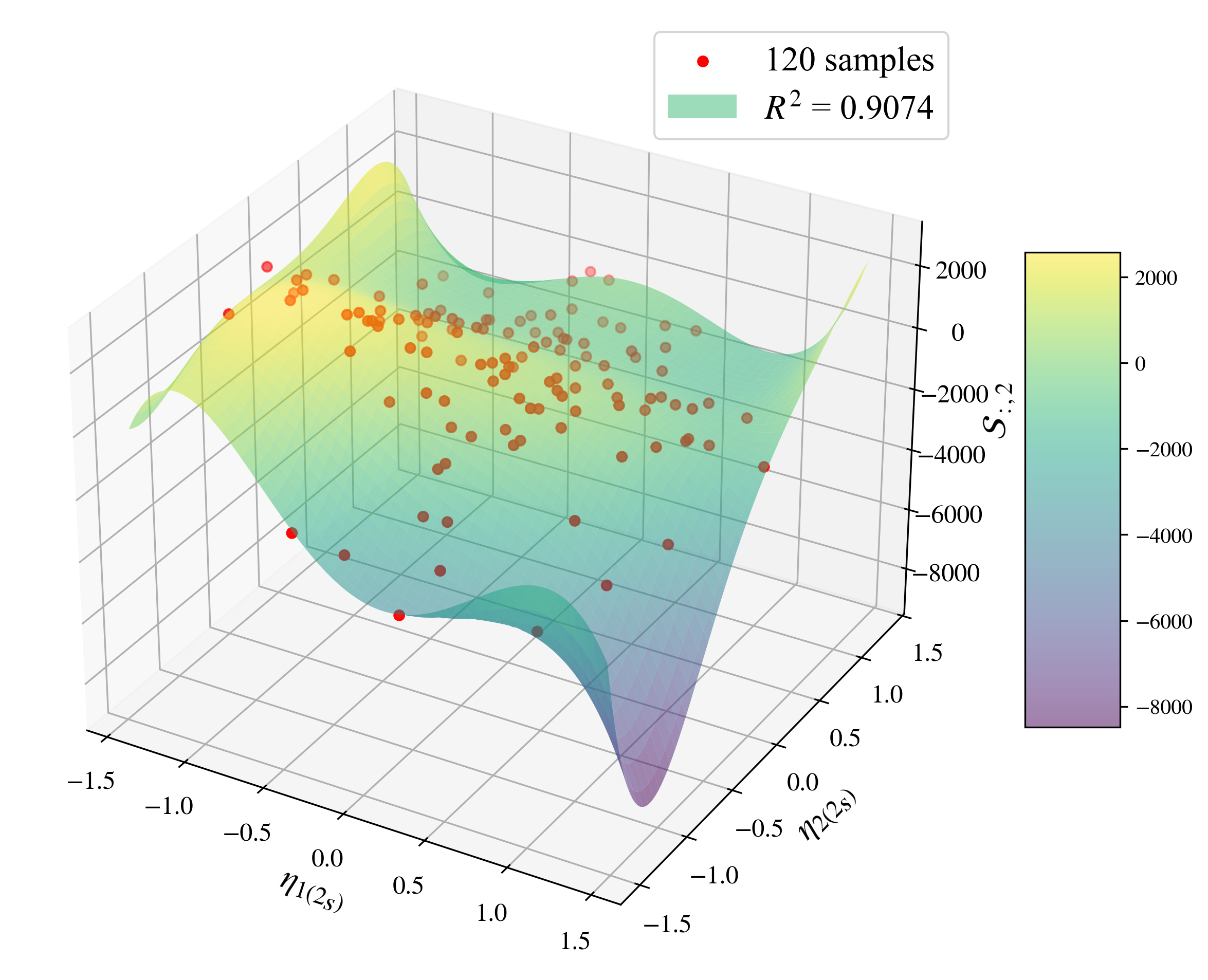}}
        \caption{Surrogate models for the first two features of the residual stress data.} 
        \label{fig: S_surrogate}
\end{figure}

We construct surrogate models for the top five features of the residual stress field. Feature 1 only uses one active variable, with $R^2$ values greater than 0.9, as shown in Figure~\ref{fig: S_surrogate}. We list the number of active variables used, polynomial type and $R^2$ value of the surrogate model for each feature in Table~\ref{tab: surrogate_S}.

\begin{table}[ht!]
\caption{Surrogate models for features 1 to 5 of residual stress, $\hat{G}_{\cS_{:, 1}}$ to $\hat{G}_{\cS_{:, 5}}$.} \label{tab: surrogate_S}
\centering
    \begin{tabular}{c c c c}
        \hline
        \textbf{Surrogate} & \textbf{Number of active variables} & \textbf{Polynomial type} & \textbf{$R^2$ Value} \\
        \hline
        $\hat{G}_{\cS_{:, 1}}$ & 1 & Parabola & 0.9665 \\
        $\hat{G}_{\cS_{:, 2}}$ & 2 & fifth-degree polynomial surface & 0.9074 \\
        $\hat{G}_{\cS_{:, 3}}$ & 2 & fifth-degree polynomial surface & 0.8869 \\
        $\hat{G}_{\cS_{:, 4}}$ & 3 & sixth-degree polynomial surface & 0.8118 \\
        $\hat{G}_{\cS_{:, 5}}$ & 2 & fifth-degree polynomial surface & 0.8493 \\
        \hline
    \end{tabular}
\end{table}

\subsection{Optimization results} \label{section: nm_opt_results}

We perform optimization with multiple initial designs using both SLSQP and COBYLA solvers in SciPy \cite{2020SciPy-NMeth} v1.1.2.  We use $\cN=20,000$ samples to estimate expectations in Algorithm~\ref{alg: solution} and set the target reliability level $\alpha_T = 0.95$. 

\begin{figure}
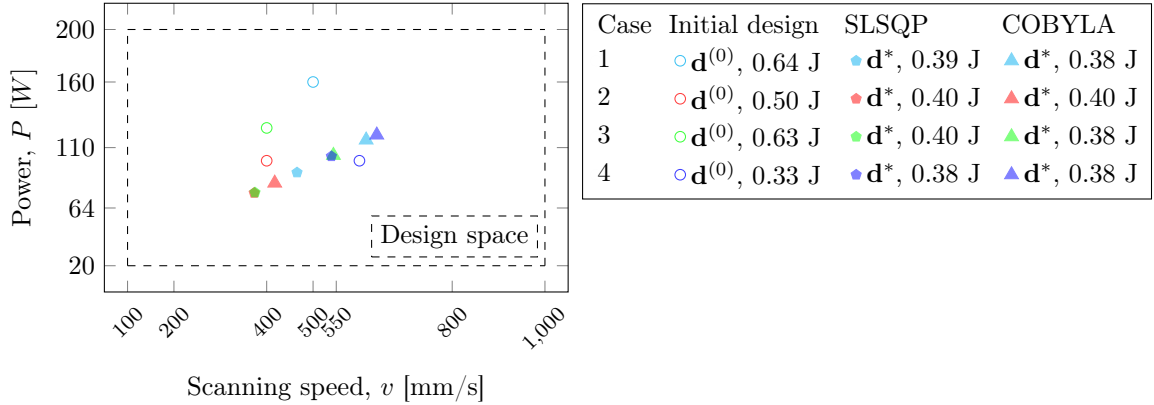

    \centering
    \include{fig_opt}
    \caption{Comparison of optimization results with different initial designs and optimization solvers. Initial design points (circles) and final optimized solutions from the SLSQP solver (pentagons) and COBYLA solver (triangles) are displayed within the bounded design space (indicated by the dashed box). All optimal $[v,P]$ solutions converge to approximately 0.4~J of consumed energy, independent of the initial guess.}
    \label{fig: opt}
\end{figure}

We start the optimization with four different initial designs. For each case, we list the initial design and the corresponding optimization result in Table~\ref{tab: optimization_results} and plot them in Figure~\ref{fig: opt}. The consumed energy is calculated using the formula $E = P \cdot t = P \cdot l/v$. In Figure~\ref{fig: opt}, we use the same color to mark each case; the initial design is marked with a circle, the optimal designs obtained by SLSQP and COBYLA solvers are marked with a pentagon and a triangle, respectively.

\begin{table}[htbp]
    \centering
    \caption{Optimization results of different initial guesses. All solutions converge to an energy consumption of approximately 0.4~J.}
    \label{tab: optimization_results}
    \begin{tabular}{c|cc|c|c|cc|ccc}
	\hline
	\multirow{3}{*}{Case} & \multicolumn{3}{c|}{Initial Designs} & \multicolumn{5}{c}{Optimization Results} \\
	\cline{2-9}
  	& \multicolumn{2}{c|}{$\bd^{\{0\}}$} & \multirow{2}{*}{$E$ [J]} 
  	& \multirow{2}{*}{Algorithm} 
  	& \multicolumn{2}{c|}{$\bd^*$} & \multirow{2}{*}{$E^*$ [J]} & \multirow{2}{*}{$\zeta^*$ [MPa]} \\
	\cline{2-3} \cline{6-7}
  	& $v$ [mm/s] & $P$ [W] & & & $v^*$ [mm/s] & $P^*$ [W] & & \\
	\hline
	\multirow{2}{*}{1} 
  	& \multirow{2}{*}{500} & \multirow{2}{*}{160} & \multirow{2}{*}{0.64} 
  	& SLSQP & 465.4 & 91.0 & 0.39 & 717.3 \\
  	& & & & COBYLA & 614.5 & 115.8 & 0.38 & 743.8 \\
\hline
\multirow{2}{*}{2} 
  	& \multirow{2}{*}{400} & \multirow{2}{*}{100} & \multirow{2}{*}{0.50} 
 	& SLSQP & 373.0 & 75.5 & 0.40 & 699.1 \\
  	& & & & COBYLA & 417.3 & 83.0 & 0.40 & 709.5 \\
\hline
\multirow{2}{*}{3} 
  	& \multirow{2}{*}{400} & \multirow{2}{*}{125} & \multirow{2}{*}{0.625} 
  	& SLSQP & 374.8 & 75.8 & 0.40 & 699.7 \\
  	& & & & COBYLA & 543.9 & 140.1 & 0.38 & 728.8 \\
	\hline
	\multirow{2}{*}{4} 
  	& \multirow{2}{*}{600} & \multirow{2}{*}{100} & \multirow{2}{*}{0.33} 
  	& SLSQP & 539.5 & 103.4 & 0.38 & 728.3 \\
  	& & & & COBYLA & 637.4 & 119.6 & 0.38 & 750.3 \\
	\hline
	\end{tabular}
\end{table}

The initial designs in case 2 and case 3 have the same scanning speed values but different beam power values. The SLSQP-optimized designs in both cases converge to nearly identical values for both parameters. In contrast, cases 2 and 4, which start with the same beam power but different scanning speeds, yield distinctly different optimal designs. This demonstrates the sensitivity of the optimization process to initial scanning speed values while suggesting robustness with respect to initial beam power selection.

In the optimization problem~\eqref{eq: optimization_reformulated}, our primary objective is to minimize energy consumption. Our analysis reveals excellent convergence behavior regarding energy consumption, despite obtaining different optimal configurations for scanning speed and beam power across all cases. The set of process parameters in the initial design of case~4 yields an energy consumption of 0.33~J. Although it is smaller than the optimal result of 0.38~J, with this parameter setting, the maximum temperature in the powder bed in the manufacturing process is $T_{\max} = 1{,}556^{\circ}$C $< T_{\text{liq}} = 1{,}650^{\circ}$C, which does not satisfy the temperature constraint in~\eqref{eq: optimization_reformulated} and indicates no powder melting. The $T_{\max}$ values corresponding to the optimal parameters obtained with SLSQP and COBYLA solvers in case~4 are both $1{,}655^{\circ}$C, which satisfies the temperature constraint. For all other cases, we see a reduction in the consumed energy from the initial designs to the optimal results. While different initial design yield varying optimal results, all solutions converge to an energy consumption of approximately 0.4~J. This consistency across different initial designs shows the robustness of our optimization approach while providing valuable guidance for selecting appropriate scanning speed and beam power.

\subsection{Validation} \label{section: nm_val}

We validate the optimal designs via their corresponding risk-based constraints on the maximum residual stress $\sigma_{\max}$. Specifically, we examine the 0.95-superquantiles of $\sigma_{\max}$ corresponding to optimal designs. The optimal auxiliary variable $\zeta^*$ can be viewed as the 0.95-quantile of $\sigma_{\max}$ related to its corresponding optimal parameter vector $\bd^*$, as discussed in Section~\ref{section: risk_measures}. We then can use $\bd^*$, $\zeta^*$, and samples of the random variable $Z$ to estimate the 0.95-superquantile with~\eqref{eq: superquantile_definition}. 

Due to computational limitations, we validate two optimization results listed in Table~\ref{tab: optimization_results}. For each optimum, we conduct 50 additional FEM simulations and compare the 0.95-superquantile from these simulations, $\bar{Q}_{0.95}^{\textsf{FEM}}$, with the 0.95-superquantile estimated from $20{,}000$ surrogate model evaluations, $\bar{Q}_{0.95}^{\textsf{Surr}}$. We generate random samples of $Z$ to create realizations of the vector $[\bd^{*\top}, Z^{\top}; \boldsymbol{\theta}^{\top}]^{\top}$, which serve as inputs for the FEM simulations and the surrogate models $\hat{G}_{\cS_j}(\boldsymbol{\eta}_{j(\bS)})$ for $j(\bS) = 1,2,3,4,5$. For FEM simulations, we obtain the residual stress field directly from the output, while for surrogate evaluations, we estimate it using~\eqref{eq: stress_estimates}. For each optimum $\bd^*$, with the samples of $\sigma_{\max}$ and the optimal auxiliary variable $\zeta^*$ corresponding to $\bd^*$, we calculate $\bar{Q}_{0.95}^{\textsf{FEM}}$ and $\bar{Q}_{0.95}^{\textsf{Surr}}$ using~\eqref{eq: superquantile_definition} and present the results in Table~\ref{tab: validation_results}. Although both validation sets only use 50 Monte Carlo samples for FEM simulations, which introduces sampling error, both superquantiles $\bar{Q}_{0.95}^{\textsf{FEM}}$ and $\bar{Q}_{0.95}^{\textsf{Surr}}$ closely match each other for both optimal results. This suggests that our approach maintains accuracy despite the limited sampling and confirms that integrating surrogate models into the optimization process produces reliable outcomes.

\begin{table}[htbp]
\centering
\caption{Validation results with the high-fidelity FEM model at two optimal solutions.}
\label{tab: validation_results}
\begin{tabular}{c|c|c|c}
\hline
\rule{0pt}{11pt} $\bd^*$ & $\bar{Q}_{0.95}^{\textsf{FEM}}$ & $\bar{Q}_{0.95}^{\textsf{Surr}}$ & $(\bar{Q}_{0.95}^{\textsf{FEM}} - \bar{Q}_{0.95}^{\textsf{Surr}})/\bar{Q}_{0.95}^{\textsf{FEM}}$ \\
\hline
$[417.3\ \text{mm/s}, 82.97\ \text{W}]$ & $768.12$ MPa & $ 762.08$ MPa & 0.79\% \\
\hline
$[614.5\ \text{mm/s}, 115.8\ \text{W}]$ & $790.69$ MPa & $804.89$ MPa & -1.80\% \\
\hline
\end{tabular}
\label{tab:validation_results}
\end{table}


\section{Conclusion and future work} \label{section: conclusion}

We formulated and solved a buffered probability of failure-constrained design optimization problem for a powder bed fusion metal additive manufacturing problem with the goal to minimize energy consumption while enhancing manufacturing quality. We simulated the manufacturing process to capture the thermal and mechanical behaviors of the manufactured part by conducting high-fidelity finite element analyses. Based on the simulation data, we developed efficient surrogate models using singular value decomposition and active subspace discovery for temporal and spatial quantities of interest, namely the temperature history and residual stress field, to reduce computational costs while maintaining prediction accuracy. Finally, we validated the optimization results through additional finite element simulations. The validation results confirm decreased energy consumption and reduced build failures.

In the future, this work can be refined with an improved two-way-coupled thermal-mechanical finite element model and a finer mesh. We furthermore plan to extend the PBOF-constrained optimization to more complex geometries by incorporating the printing path into the optimization formulation.

\section*{Statements and declarations}

\subsection*{Acknowledgment}
Funding information: Y. Guo and B. Kramer were financially supported by the Air Force Office of Scientific Research (PM Fahroo) under award number FA9550-24-1-0237 and the Defense Advanced Research Projects Agency (DARPA) Cooperative Agreement No. HR0011-25-2-0009, "Predictive Real-time Intelligence for Metallic Endurance (PRIME)."

\subsection*{CRediT authorship contribution statement}
Yulin Guo: Conceptualization, Methodology, Data Curation, Software, Formal analysis, Investigation, Software, Visualization, Writing - original draft; Boris Kramer: Conceptualization, Formal analysis, Writing - review and editing, Funding acquisition, Project administration, Supervision. 
\subsection*{Conflict of interest}
Boris Kramer reports a relationship with ASML Holding US that includes: consulting or advisory. Yulin Guo declares that he has no known competing financial interests or personal relationships that could have appeared to influence the work reported in this paper.
\subsection*{Data availability}
Input files for \textsf{Abaqus} simulations and the complete data sets used to generate the results are available to download from the Github repository: https://github.com/yulin-g/BPOF-optimization-AM. 
\subsection*{Replication of results}
The results presented in this manuscript can be fully replicated using the provided material in the Github repository.
\subsection*{Ethics approval and consent to participate}
Not applicable.
\bibliography{references}
\bibliographystyle{acm}

\end{document}

%% file: mymacros.tex
\newcommand{\bit}{\begin{itemize}}
\newcommand{\eit}{\end{itemize}}
\newcommand{\ben}{\begin{enumerate}}
\newcommand{\een}{\end{enumerate}}

\newcommand {\real} {\mathbb{R}}

\newcommand{\bF}{\ensuremath{\mathbf{F}}}

\newcommand{\bS}{\ensuremath{\mathbf{S}}}
\newcommand{\bT}{\ensuremath{\mathbf{T}}}
\newcommand{\bU}{\ensuremath{\mathbf{U}}}
\newcommand{\bV}{\ensuremath{\mathbf{V}}}
\newcommand{\bW}{\ensuremath{\mathbf{W}}}

\newcommand{\bd}{\ensuremath{\mathbf{d}}}

\renewcommand{\bf}{\ensuremath{\mathbf{f}}}

\newcommand{\bu}{\ensuremath{\mathbf{u}}}
\newcommand{\bv}{\ensuremath{\mathbf{v}}}
\newcommand{\bw}{\ensuremath{\mathbf{w}}}

\newcommand{\cC}{\ensuremath{\mathcal{C}}}
\newcommand{\cD}{\ensuremath{\mathcal{D}}}
\newcommand{\cF}{\ensuremath{\mathcal{F}}}

\newcommand{\cM}{\ensuremath{\mathcal{M}}}
\newcommand{\cN}{\ensuremath{\mathcal{N}}}

\newcommand{\cS}{\ensuremath{\mathcal{S}}}
\newcommand{\cT}{\ensuremath{\mathcal{T}}}



%% file: fig_Q_alpha_POF_all.tex
\pgfmathdeclarefunction{frechet}{3}{%
	\pgfmathparse{(#1/#2)*((x-#3)/#2)^(-1-#1)*exp(-((x-#3)/#2)^(-#1))}%
}
\newcommand{\pkg}[1]{\texttt{#1}}
\newcommand{\cls}[1]{\textsf{#1}}
\newcommand{\file}[1]{\texttt{#1}}

\newcommand{\failset}{\mathcal{G}}

\newcommand{\qbar}{\overline{Q}}
\newcommand{\pt}{p_t}
\newcommand{\hpt}{\hat{p}_t}
\newcommand{\pbuf}{\overline{p}_t}
\begin{tikzpicture}
\begin{groupplot}[
    group style={
        group size = 3 by 1,
        horizontal sep = 1cm,
    },
    width=5cm, height=3.2cm,
    xmin=0.75, xmax=3.8, ymin=0, ymax=0.7,
    samples=100,
    axis lines=left,
    every axis x label/.style={at=(current axis.right of origin),anchor=north},
    xlabel={$g(\bd,Z;\boldsymbol{\theta})$},
    every tick/.style={black, very thick},
    tick label style={font=\footnotesize,black,text height=1ex},
    ytick=\empty,
    enlargelimits=upper,
    clip=true,
    axis on top
]

\nextgroupplot[
    ylabel={Probability\\density},
    every axis y label/.style={
        anchor=south west,
        rotate=90,
        align=center,
        inner sep=0.4ex,
        xshift=0cm
    },
    xtick={2.3,2.6},
    xticklabels={$Q_{\alpha_T}$,\hspace*{2.5ex}$\tau$}
]
\addplot [fill=blue!30, draw=none, domain=2.6:3.5] {frechet(3,2,0)} \closedcycle;
\addplot [very thick, black!80, smooth, domain=0.3:3.5] {frechet(3,2,0)};    
\node at (axis cs:3.2,0.5) [fill=white,font=\footnotesize,inner sep=0.18ex] {$p_{\tau} < 1-\alpha_T$};

\nextgroupplot[
    xtick={2.4},
    xticklabels={$Q_{\alpha_T}= \tau$}
]
\addplot [fill=blue!30, draw=none, domain=2.4:3.5] {frechet(3,2,0)} \closedcycle;
\addplot [very thick, black!80, smooth, domain=0.3:3.5] {frechet(3,2,0)};    
\node at (axis cs:3.2,0.5) [fill=white,font=\footnotesize,inner sep=0.18ex] {$p_{\tau} = 1-\alpha_T$};

\nextgroupplot[
    xtick={2.3,2.6},
    xticklabels={$\tau$, \hspace*{2.5ex} $Q_{\alpha_T}$}
]
\addplot [fill=blue!30, draw=none, domain=2.3:3.5] {frechet(3,2,0)} \closedcycle;
\addplot [very thick, black!80, smooth, domain=0.3:3.5] {frechet(3,2,0)};    
\node at (axis cs:3.2,0.5) [fill=white,font=\footnotesize,inner sep=0.18ex] {$p_{\tau} > 1-\alpha_T$};

\end{groupplot}
\end{tikzpicture}

%% file: fig_BPOF_buffer.tex
\pgfmathdeclarefunction{frechet}{3}{%
	\pgfmathparse{(#1/#2)*((x-#3)/#2)^(-1-#1)*exp(-((x-#3)/#2)^(-#1))}%
}
\newcommand{\pkg}[1]{\texttt{#1}}
\newcommand{\cls}[1]{\textsf{#1}}
\newcommand{\file}[1]{\texttt{#1}}

\newcommand{\failset}{\mathcal{G}}

\newcommand{\qbar}{\overline{Q}}
\newcommand{\pt}{p_{\tau}}
\newcommand{\hpt}{\hat{p}_{\tau}}
\newcommand{\pbuf}{\overline{p}_{\tau}}
\begin{tikzpicture}[scale=1,auto]
\begin{axis}[
    xmin=0.75, xmax=5.4, ymin=0, ymax=0.7, height=5cm, width=8cm,
    black!60, very thick, domain=0:2.5, samples=100,
    axis lines=left,
    every axis y label/.style={
        at={(ticklabel cs:0.5)},
        rotate=90,
        align=center, 
        yshift = 5pt
        },
    ylabel={Probability density},
    every axis x label/.style={at=(current axis.right of origin),anchor=north}, 
    xlabel={$g(\bd,Z;\boldsymbol{\theta})$}, 
    every tick/.style={black, very thick}, 
    tick label style={font=\footnotesize,black,text height=1ex}, 
    xtick={2.4,3.5}, xticklabels={$Q_{\alpha} = \zeta$,\hspace*{2.5ex}$\bar{Q}_{\alpha} = \tau$}, 
    ytick=\empty, enlargelimits=upper, clip=true, axis on top, smooth
]
\addplot [fill=blue!30, draw=none, domain=2.4: 3.5] {frechet(3,2,0)} \closedcycle;
\addplot [fill=red!30, draw=none, domain=3.5:5] {frechet(3,2,0)} \closedcycle;
\addplot [very thick,black!80, domain=0.3:5] {frechet(3,2,0)};	

\node at (axis cs:3.45,0.44) [fill=white, font=\footnotesize, inner sep=0.18ex, text=blue, opacity=0.6] {$\mathbb{P}[g(\bd,Z;\boldsymbol{\theta}) \in [\zeta, \tau]]$};
\node at (axis cs:4.6,0.25) [fill=white, font=\footnotesize, inner sep=0.18ex, text=red, opacity=0.6] {$p_{\tau} = \mathbb{P}[g(\bd,Z;\boldsymbol{\theta}) > \tau]$};
\node at (axis cs:2.2,0.65) [fill=white, font=\footnotesize, inner sep=0.18ex] {$\bar{p}_{\tau} = $};
\node at (axis cs:3.45,0.65) [fill=white, font=\footnotesize, inner sep=0.18ex, text=blue, opacity=0.6] {$\mathbb{P}[g(\bd,Z;\boldsymbol{\theta}) \in [\zeta, \tau]]$};
\node at (axis cs:4.55,0.65) [fill=white, font=\footnotesize, inner sep=0.18ex] {+};
\node at (axis cs:4.8,0.64) [fill=white, font=\footnotesize, inner sep=0.18ex, text=red, opacity=0.6] {$p_{\tau}$};
\node at (axis cs:5.35,0.65) [fill=white, font=\footnotesize, inner sep=0.18ex] {$ = 1 - \alpha$};
    
    
    
    
    
    
    \node[draw, fill=white, align=center] at (2.9,0.08) {buffer};
    
    \end{axis}	
\end{tikzpicture}

%% file: fig_T_L_2.tex
\begin{tikzpicture}
\begin{axis}[
    width=\linewidth,
    height=0.6\linewidth,
    grid=both,
    grid style={dashed,gray!30},
    xlabel={Number of singular value-singular vector pairs, $k$},
    ymode = log,
    ylabel={$\textsf{err}_{\bT}(k)$},
    ymode = log,
    scaled y ticks=false,
    ytick={0.1, 0.01, 0.001},         
    yticklabels={$10^{-1}$, $10^{-2}$, $10^{-3}$},  
    ymajorgrids=true,                 
]
]

\addplot[
    color=blue,
    mark=*,
    mark size=2pt,
    thick
] coordinates {
    (1, 0.10325313944291051)
    (2, 0.04530458923079103)
    (3, 0.008644761206775068)
    (4, 0.0032871406645308397)
    (5, 0.001943224380352817)
    (6, 0.0015067953784717004)
    (7, 0.0012720758890423156)
    (8, 0.0010631241969934668)
    (9, 0.0008968229291920379)
    (10, 0.0007490274357354919)
};

\addplot[
    color=blue,
    dashed,
    thick
] coordinates {
    (1, 0.10325313944291051)
    (2, 0.04530458923079103)
    (3, 0.008644761206775068)
    (4, 0.0032871406645308397)
    (5, 0.001943224380352817)
    (6, 0.0015067953784717004)
    (7, 0.0012720758890423156)
    (8, 0.0010631241969934668)
    (9, 0.0008968229291920379)
    (10, 0.0007490274357354919)
};


\node[circle, fill=lime, inner sep=3pt] (pt) at (axis cs:2, 0.0453) {};

\node[align=center, fill=white, anchor=west] (label1) at (axis cs:3, 0.04530458923079103) {$\textsf{err}_{\bT}(2) = 4.53\% < 5\%$};

\node[align=center, fill=white, anchor=north] (label2) at (label1.south) {$K_{\cT} = 2$};

\draw[<-] (pt) -- (label1);

\end{axis}
\end{tikzpicture}

%% file: fig_S_L_2.tex
\begin{tikzpicture}
\begin{axis}[
    width=\linewidth,
    height=0.6\linewidth,
    grid=both,
    grid style={dashed,gray!30},
    xlabel={$k$},
    ymode = log,
    ylabel={$\textsf{err}_{\bS}(k)$},
    scaled y ticks=false,
    ytick={1, 0.1},         
    yticklabels={$10^{0}$, $10^{-1}$},  
    ymajorgrids=true,                 
    ymax = 1.5
]

\addplot[
    color=blue,
    mark=*,
    mark size=2pt,
    thick
] coordinates {
    (1, 0.39987282961628545)
    (2, 0.26315456388152675)
    (3, 0.15059725663747497)
    (4, 0.13418503344473287)
    (5, 0.10370120666329985)
    (6, 0.09284355897732577)
    (7, 0.08453758821970234)
    (8, 0.07918093976700538)
    (9, 0.07501744014285021)
    (10, 0.07074871407171848)
};

\addplot[
    color=blue,
    dashed,
    thick
] coordinates {
    (1, 0.39987282961628545)
    (2, 0.26315456388152675)
    (3, 0.15059725663747497)
    (4, 0.13418503344473287)
    (5, 0.10370120666329985)
    (6, 0.09284355897732577)
    (7, 0.08453758821970234)
    (8, 0.07918093976700538)
    (9, 0.07501744014285021)
    (10, 0.07074871407171848)
};


\node[circle, fill=lime, inner sep=3pt] (pt) at (axis cs:5, 0.10370120666329985) {};

\node[circle, fill=lime, inner sep=1pt] (pt2) at (axis cs:5, 0.22) {};

\node[align=center, fill=white] (label1) at (axis cs:5.8, 0.28) {$K_{\cS} = 5$};

\node[align=center, fill=white, anchor=south] (label2) at (label1.north) {$\textsf{err}_{\bS}(5)$ = 10.37\% $\approx$ 10\%};

\node[align=left, fill=white, anchor=west] at (axis cs:6, 0.15) {$\textsf{err}_{\bS}(6)$ = 9.28\%};

\draw[<-] (pt) -- (pt2);

\end{axis}
\end{tikzpicture}

%% file: fig_T_surrogate_1.tex
\begin{tikzpicture}
\begin{axis}[
    width=\linewidth,
    height=0.8\linewidth,
    xlabel={$\eta_{1(1_{\bT})}$},
    ylabel={$\cT_{:, 1}$},
    grid=both,
    minor grid style={gray!25},
    major grid style={gray!25},
    xmin=-1.5, xmax=1.3,
    ymin=-12500, ymax=-3500,
    xtick={-1.5, -1.0, -0.5, 0.0, 0.5, 1.0},
    ytick={-12000, -10000, -8000, -6000, -4000},
    tick scale binop=\times,
    legend pos=south west,
    legend style={draw=none, fill=white, fill opacity=0.8}
]

\addplot[only marks, mark=*, mark size=1.5pt, color=red] coordinates {
(-0.35604601, -5222.68146384)
(-0.85321957, -4431.24858071)
(-0.15003523, -6366.75728066)
(0.64428471, -8488.00419495)
(-0.3573726, -5778.08280897)
(-0.01682151, -6517.37002178)
(0.30224563, -7042.31621676)
(-0.6375322, -4509.5260367)
(-0.16125536, -5932.50714541)
(0.1671072, -7173.47868631)
(1.21382304, -11985.25341959)
(-0.4594446, -5525.208337)
(0.55384006, -8501.70912273)
(-0.78402307, -4872.07050438)
(0.21205757, -7726.73598906)
(-1.06530757, -4396.21306846)
(-1.10836451, -4322.53074399)
(-0.75421247, -5336.37819937)
(-0.98384585, -4118.66112941)
(-0.02408545, -6363.82222159)
(-0.24540995, -5435.19348756)
(1.19080491, -11501.10026209)
(-0.27084559, -5032.147315)
(-0.18641305, -6343.28003547)
(0.43598611, -8003.54591547)
(0.12280081, -7053.66055199)
(0.09594662, -6921.99197648)
(0.10180872, -7049.17199493)
(0.75292793, -8985.2159661)
(-0.66019518, -5323.23195922)
(0.15481499, -6969.72313903)
(-0.31610118, -6004.01320473)
(0.38159276, -7436.31492562)
(0.35355221, -7585.87753021)
(0.16245266, -7271.72128137)
(0.32958283, -7876.07792871)
(-0.00958807, -6603.04571661)
(0.46337261, -8403.67145474)
(-0.26870311, -6187.00619045)
(-0.96455013, -4350.96064677)
(0.13601139, -6988.34973093)
(0.66687046, -8914.41252307)
(0.95314675, -11389.13590645)
(0.54707974, -7937.95792349)
(0.00397185, -6665.43926937)
(-0.51812332, -5663.69932116)
(0.52154515, -7972.61927832)
(-0.13316052, -6434.27307711)
(-0.60487159, -5370.58182302)
(0.73618303, -8954.04782815)
(0.09214349, -7165.677098)
(-0.1787798, -6025.81010056)
(-0.19801392, -5654.56207964)
(0.57064255, -8624.03480113)
(0.92455955, -9407.21779074)
(0.41982038, -8085.84259841)
(0.10706527, -6870.83406993)
(-0.39263227, -6089.02584051)
(0.31475638, -7441.5970808)
(-0.36269937, -4919.48961432)
(-0.38670171, -5165.92217136)
(0.01297296, -6592.21850529)
(-0.71152114, -4820.02649833)
(0.84557608, -9283.03487533)
(-0.32794099, -5901.38363433)
(0.19969281, -7215.84528957)
(-0.26150438, -5988.92215278)
(-1.26027646, -3977.67156267)
(1.0832421, -10341.51789066)
(-0.67600595, -4553.23930083)
(-0.46434135, -5462.89746507)
(0.89940882, -9289.98764959)
(0.17192763, -7228.27805222)
(0.73412719, -10434.60232508)
(-0.25259304, -6384.21593879)
(-1.32404301, -4053.49750739)
(-0.84486431, -4958.90743362)
(0.24268422, -7452.29890768)
(-0.43049652, -5788.30106194)
(-0.31996468, -5633.87714562)
(0.15514527, -7229.23287569)
(0.13998067, -6955.06272774)
(0.96672985, -10313.06109206)
(0.97071801, -10604.61103261)
(0.6771973, -9396.83708829)
(-0.4262096, -5872.4419686)
(0.6680305, -10340.92638003)
(0.98492091, -10252.94559178)
(0.08610858, -6835.58226899)
(-0.13737846, -6084.10350195)
(-0.29922473, -5962.51870261)
(-0.20448, -6386.99498652)
(0.14051913, -6935.75429334)
(0.51132955, -9254.19486296)
(0.45763613, -8192.25524548)
(0.08462817, -6929.03396432)
(-0.75448754, -5084.27856246)
(0.21915322, -7025.99803879)
(0.14221677, -7294.34865133)
(-0.61897004, -5573.92928382)
(-0.12834356, -6338.24193153)
(-0.00921877, -6424.97344975)
(-0.3031132, -5684.23728475)
(0.7868163, -9169.80248251)
(-0.37113817, -5649.72433225)
(-0.70173016, -4758.27178951)
(-0.82804717, -4276.8957791)
(-0.57046729, -5577.29190452)
(0.49303287, -7727.27374873)
(-0.44604354, -5249.20194323)
(-0.17407323, -5909.10948884)
(0.39273361, -7562.93041139)
(-0.9156727, -4074.08693612)
(-0.28827763, -5990.72435102)
(0.53285977, -8382.60378382)
(0.53164719, -7947.56516146)
(0.09443822, -6983.46594001)
(-0.01456685, -6592.55103151)
(0.77771827, -9685.0785172)
(-0.22781608, -5662.29381938)
};
\addlegendentry{120 samples}

\addplot[dashed, thick, color=red, domain=-1.5:1.3] {-6870.126942468027 + x * (-3143.7)};

\addlegendentry{Surrogate, $R^2 = 0.9245$}

\end{axis}
\end{tikzpicture}

%% file: fig_T_surrogate_2.tex
\begin{tikzpicture}
\begin{axis}[
    width=\linewidth,
    height=0.8\linewidth,
    xlabel={$\eta_{1(2_{\bT})}$},
    ylabel={$\cT_{:, 2}$},
    grid=both,
    minor grid style={gray!25},
    major grid style={gray!25},
    xmin=-1.6, xmax=1.5,
    ymin=-2400, ymax=1200,
    xtick={-1.5, -1.0, -0.5, 0.0, 0.5, 1.0, 1.5},
    ytick={-2000, -1000, 0, 1000},
    tick scale binop=\times,
    legend pos=south west,
    legend style={draw=none, fill=white, fill opacity=0.8}
]

\addplot[only marks, mark=*, mark size=1.5pt, color=red] coordinates {
(-0.09574095, 565.76663826)
(-0.73306763, 791.45644601)
(-0.46126641, 650.67872563)
(0.53912374, -209.08565076)
(-0.37918656, 638.28074998)
(0.3374514, -53.52609186)
(-0.1348934, 690.04750971)
(-0.3773338, 734.73749696)
(0.08609603, 301.46426744)
(0.16766661, 161.90487632)
(1.22338111, -2083.47548565)
(-0.45648478, 668.91981239)
(0.567212, -358.60086111)
(-0.75294054, 790.4910733)
(0.56390421, -564.15881903)
(-1.04554961, 812.79626493)
(-1.1275584, 917.79646693)
(-0.91493143, 718.74872142)
(-0.85115523, 877.18133689)
(0.31067185, 48.88556461)
(0.12844556, 365.82934413)
(1.12142307, -1687.41107421)
(0.19106826, 406.62137141)
(-0.48494181, 675.22671666)
(0.45180055, -178.0300249)
(0.10174756, 222.0007439)
(-0.06973608, 377.08597292)
(0.06366012, 169.31826264)
(0.59973069, -393.28307419)
(-0.74516281, 696.55507119)
(-0.25007997, 633.10150728)
(-0.36891086, 520.48438453)
(-0.01918191, 465.37190148)
(-0.00290425, 225.82400758)
(-0.05082577, 226.48258144)
(0.48005635, -293.72631182)
(0.34073193, -89.85030728)
(0.62057796, -558.184899)
(-0.58517129, 724.26163323)
(-0.85671299, 789.40783147)
(-0.06584904, 406.43118671)
(0.66277732, -500.53950315)
(1.02273648, -1907.97067489)
(0.21880928, 213.41488655)
(0.06725452, 254.39453985)
(-0.68972796, 743.30613482)
(0.30800293, 118.9918825)
(-0.36148903, 565.20260785)
(-0.72337913, 804.04697285)
(0.6301831, -441.21505631)
(0.44189497, -341.98987528)
(-0.07375244, 462.26309865)
(0.10040208, 351.68026535)
(0.54323064, -451.43004503)
(0.71724852, -514.4613023)
(0.45366036, -234.77280405)
(0.01199014, 359.9083839)
(-0.6432316, 572.50044951)
(0.1985887, 192.77446563)
(-0.00263791, 575.27459033)
(-0.10156216, 524.32992879)
(0.06539166, 352.89197941)
(-0.5824612, 700.81270381)
(0.72575179, -554.96489467)
(-0.42540277, 659.97081067)
(0.27772509, 83.52216408)
(-0.17181829, 385.55402239)
(-1.20871577, 886.40679921)
(1.02888313, -1059.04670832)
(-0.46987408, 786.03437306)
(-0.41998668, 646.0568753)
(0.70673745, -454.95387321)
(-0.24216287, 377.40443778)
(0.94406517, -1650.77780144)
(-0.45113576, 470.9395432)
(-1.34521181, 906.36221839)
(-0.91380567, 768.10362481)
(0.29647705, 4.2890297)
(-0.51886162, 585.77850062)
(-0.16653857, 516.86420882)
(0.37540592, -158.68035369)
(-0.14731569, 524.03601861)
(0.97863335, -1172.01213733)
(0.99923787, -1384.24796481)
(0.74950641, -865.81249253)
(-0.55862754, 586.96759808)
(0.88172757, -1620.03124875)
(0.95300508, -1120.66654565)
(0.25668118, 123.43029079)
(-0.01479368, 467.85912564)
(-0.36935935, 604.67879916)
(-0.44422519, 567.91954139)
(-0.10832244, 467.60164461)
(0.79368809, -1221.73767452)
(0.47116455, -245.45699511)
(0.00293951, 266.9293236)
(-0.83935322, 788.18180172)
(-0.01517367, 542.53660777)
(0.49133538, -368.64030683)
(-0.83825528, 699.04282776)
(-0.50621101, 774.5289033)
(0.36879412, -80.43963335)
(-0.16853136, 538.9393418)
(0.6372023, -538.25958062)
(-0.30568255, 572.09129737)
(-0.56172601, 718.79840802)
(-0.62357562, 781.66295266)
(-0.70798761, 684.77722842)
(0.21646782, 273.84585314)
(-0.24674434, 572.94446982)
(0.10510568, 285.64114311)
(0.21550716, 241.7073333)
(-0.73866396, 870.07422343)
(-0.27762897, 459.7664548)
(0.56618522, -305.15873805)
(0.36313187, 109.28808682)
(-0.19345641, 422.10874159)
(0.27847587, 9.96131042)
(0.81031965, -964.46391004)
(0.11980038, 315.49586943)

};
\addlegendentry{120 samples}

\addplot[dashed, thick, color=red, domain=-1.5:1.3] {-662.9651309 * x^2 - 1063.48001925 * x + 368.67147444};

\addlegendentry{Surrogate, $R^2 = 0.9836$}

\end{axis}
\end{tikzpicture}

%% file: fig_S_surrogate_1.tex
\begin{tikzpicture}
\begin{axis}[
    width=\linewidth,
    height=0.8\linewidth,
    xlabel={$\eta_{1(1_{\bS})}$},
    ylabel={$\cS_{:, 1}$},
    grid=both,
    minor grid style={gray!25},
    major grid style={gray!25},
    xmin=-1.6, xmax=1.5,
    ymin=-14500, ymax=1000,
    xtick={-1.5, -1.0, -0.5, 0.0, 0.5, 1.0, 1.5},
    tick scale binop=\times,
    legend pos=south west,
    legend style={draw=none, fill=white, fill opacity=0.8}
]

\addplot[only marks, mark=*, mark size=1.5pt, color=red] coordinates {
(-0.35407406, -710.21404731)
(-0.80928645, -1505.6399499)
(-0.25088211, -2261.37657319)
(0.465488, -8349.12348202)
(-0.47281281, -989.32187302)
(-0.15291872, -2101.82419295)
(0.2440657, -5413.53280881)
(-0.71011357, -355.38545158)
(-0.26400037, -2075.11955258)
(0.21482498, -5556.71651131)
(1.19074082, -12517.51104327)
(-0.66535757, -490.1730341)
(0.42986847, -8161.07863389)
(-0.69960229, -691.95426328)
(0.08434977, -5646.89566232)
(-0.78884714, -939.50446508)
(-1.17245886, -451.26551199)
(-0.39560517, -1700.52528515)
(-0.95743867, -104.08145676)
(0.01152173, -2666.7235764)
(-0.47924602, -936.56561609)
(1.2247028, -11194.9469184)
(-0.24493529, -2088.87959364)
(0.03713031, -4664.90653823)
(0.54869653, -8836.02748146)
(0.14841208, -5588.85893478)
(0.20100825, -5867.70981387)
(0.18923497, -5638.37573982)
(0.7976017, -9387.60725607)
(-0.64118813, -1393.59119518)
(0.07456119, -4536.59507275)
(-0.32722885, -1241.69141297)
(0.56870581, -7949.76763613)
(0.27462172, -5773.64361909)
(0.38036648, -6916.82718098)
(0.18669051, -6523.18589366)
(-0.07769902, -3143.80497354)
(0.52323257, -9346.75019021)
(0.01186641, -4194.30877762)
(-0.99545775, -691.75049869)
(0.25562773, -6493.89776557)
(0.67919809, -10059.49523044)
(0.79477108, -11422.15404854)
(0.59881747, -7878.37180248)
(-0.07996263, -2899.97224745)
(-0.44094345, -1319.41966287)
(0.43337259, -7492.31426988)
(-0.20870787, -2866.8710102)
(-0.40448459, -1713.20782481)
(0.7811846, -9802.60311818)
(-0.15400783, -3332.86497669)
(-0.21911998, -1476.13731157)
(-0.2567525, -1846.37340308)
(0.65231448, -9402.33782839)
(0.91108033, -9902.12309079)
(0.35290636, -6997.81334708)
(0.20766012, -5587.54288945)
(-0.26801274, -2573.52109491)
(0.3072119, -6855.26367269)
(-0.34528393, -1117.94855408)
(-0.4282425, -1040.83060484)
(-0.10711248, -2513.1134798)
(-0.86638, -691.24491918)
(0.83244646, -9575.99487144)
(-0.35253241, -2115.20973756)
(0.32271453, -7124.31832478)
(-0.29745857, -1149.06653465)
(-1.0415948, -800.67186138)
(1.102144, -10865.56241678)
(-0.81401493, -727.82902543)
(-0.53749621, -605.97268957)
(0.9910711, -10736.14759763)
(0.38252769, -6496.14062874)
(0.6129095, -11352.33652488)
(-0.20283206, -2695.28059559)
(-1.24638425, -734.08730806)
(-0.59049858, -1800.50832615)
(0.13076608, -5954.68986457)
(-0.50209115, -1142.21558443)
(-0.32731799, -550.45885739)
(0.18391917, -6469.313361)
(0.11709358, -4999.22701181)
(1.07160981, -10388.23846151)
(0.87276282, -10769.3746215)
(0.65939439, -9935.78886315)
(-0.52653234, -1187.19190093)
(0.6861273, -11340.22069592)
(1.13387108, -10737.65314164)
(0.04649235, -4225.6829543)
(-0.10292285, -2391.18874136)
(-0.49753788, -577.83943003)
(-0.16247714, -2956.49603051)
(0.29285433, -6308.92175523)
(0.3839049, -9302.11508471)
(0.58327746, -8851.66583798)
(0.14059922, -4853.71230491)
(-0.79099113, -1268.32854235)
(0.51827084, -7748.19983707)
(0.12875592, -4789.72890601)
(-0.68556474, -758.79935771)
(-0.13515209, -2879.46067546)
(-0.30567659, -757.20925649)
(-0.28542143, -1642.84632961)
(0.75212139, -9109.31611668)
(-0.43306524, -344.21649852)
(-0.6134659, -1336.66264265)
(-0.72608942, -64.21501921)
(-0.4810837, -1281.28135044)
(0.36257712, -6642.58222038)
(-0.47038756, -703.36415759)
(-0.44283246, -529.62986884)
(0.53721382, -8031.96486146)
(-0.81662533, -1040.42980021)
(-0.47129856, -703.08810709)
(0.53244211, -9131.49858493)
(0.36259561, -6839.43615583)
(0.29909953, -6514.57695068)
(-0.11579535, -2877.60350748)
(0.60481644, -9373.30008503)
(-0.2720456, -695.81250311)

};
\addlegendentry{120 samples}

\addplot[dashed, thick, color=red, domain=-1.5:1.3] {-1584.85882257 * x^2 -6015.99306103 * x + -4075.35568946};

\addlegendentry{Surrogate, $R^2 = 0.9665$}

\end{axis}
\end{tikzpicture}

%% file: fig_opt.tex
\newcommand{\indentedLegend}[2]{\hspace{#1}\,#2}

\begin{tikzpicture}
    \begin{axis}[
        width=0.5\linewidth,
        height=0.35\linewidth,
        xlabel={Scanning speed, $v$ [mm/s]},
        ylabel={Power, $P\ [W]$},
        xtick={100, 200, 400, 500, 550, 800, 1000},
        ytick={20, 64, 110, 160, 200},
        xticklabel style={rotate=45, font = \footnotesize},
        legend columns=4,
        legend cell align={left},
        legend style={fill opacity=0.8, draw opacity=1, text opacity=1, legend pos=outer north east},
        xmin=50,    
    	xmax=1050, 
    	ymin=0,    
    	ymax=220   
    ]
    
    \addplot[only marks, mark=o, color=cyan, opacity=0.75, mark size=2, mark options={solid, fill opacity=0.005, draw opacity=0.009}] coordinates {(1500, 250)};
    \addlegendentry{\indentedLegend{-0.5em}{Case $\ $}}
    \addplot[only marks, mark=o, color=cyan, opacity=0.75, mark size=2, mark options={solid, fill opacity=0.005, draw opacity=0.009}] coordinates {(1500, 250)};
    \addlegendentry{\indentedLegend{-1em}{Initial design $\ \ $}}
    \addplot[only marks, mark=pentagon*, color=cyan, opacity=0.5, mark size=2, mark options={solid, fill opacity=0.005, draw opacity=0.009}] coordinates {(1500, 250)};
    \addlegendentry{\indentedLegend{-1em}{SLSQP $\quad \quad  \ $}}
    \addplot[only marks, mark=triangle*, color=cyan, opacity=0.5, mark size=3, mark options={solid, fill opacity=0.005, draw opacity=0.009}] coordinates {(1500, 250)};
    \addlegendentry{\indentedLegend{-1em}{COBYLA $\ $}}
    
    \addplot[only marks, mark=o, color=cyan, opacity=0.75, mark size=2, mark options={solid, fill opacity=0.005, draw opacity=0.009}] coordinates {(1500, 250)};
    \addlegendentry{\indentedLegend{-0.5em}{1}}
    \addplot[only marks, mark=o, color=cyan, opacity=0.75, mark size=2, mark options={solid}] coordinates {(500,160)};
    \addlegendentry{$\bd^{(0)}$, 0.64 J }
    \addplot[only marks, mark=pentagon*, color=cyan, opacity=0.5, mark size=2] coordinates {(465.4,91.04)};
    \addlegendentry{$\bd^*$, 0.39 J }
    \addplot[only marks, mark=triangle*, color=cyan, opacity=0.5, mark size=3] coordinates {(614.5,115.8)};
    \addlegendentry{$\bd^*$, 0.38 J }
    
    \addplot[only marks, mark=o, color=cyan, opacity=0.75, mark size=2, mark options={solid, fill opacity=0.005, draw opacity=0.009}] coordinates {(1500, 250)};
    \addlegendentry{\indentedLegend{-0.5em}{2}}
    \addplot[only marks, mark=o, color=red, opacity=0.75, mark size=2, mark options={solid}] coordinates {(400,100)};
    \addlegendentry{$\bd^{(0)}$, 0.50 J }
    \addplot[only marks, mark=pentagon*, color=red, opacity=0.5, mark size=2] coordinates {(373,75.48)};
    \addlegendentry{$\bd^*$, 0.40 J }
    \addplot[only marks, mark=triangle*, color=red, opacity=0.5, mark size=3] coordinates {(417.3,82.97)};
    \addlegendentry{$\bd^*$, 0.40 J }

    \addplot[only marks, mark=o, color=cyan, opacity=0.75, mark size=2, mark options={solid, fill opacity=0.005, draw opacity=0.009}] coordinates {(1500, 250)};
    \addlegendentry{\indentedLegend{-0.5em}{3}}
    \addplot[only marks, mark=o, color=green, opacity=0.75, mark size=2, mark options={solid}] coordinates {(400,125)};
    \addlegendentry{$\bd^{(0)}$, 0.63 J }
    \addplot[only marks, mark=pentagon*, color=green, opacity=0.5, mark size=2] coordinates {(374.8,75.8)};
    \addlegendentry{$\bd^*$, 0.40 J }
    \addplot[only marks, mark=triangle*, color=green, opacity=0.5, mark size=3] coordinates {(543.9,104.1)};
    \addlegendentry{$\bd^*$, 0.38 J }

    \addplot[only marks, mark=o, color=cyan, opacity=0.75, mark size=2, mark options={solid, fill opacity=0.005, draw opacity=0.009}] coordinates {(1500, 250)};
    \addlegendentry{\indentedLegend{-0.5em}{4}}
    \addplot[only marks, mark=o, color=blue, opacity=0.75, mark size=2, mark options={solid}] coordinates {(600,100)};
    \addlegendentry{$\bd^{(0)}$, 0.33 J }
    \addplot[only marks, mark=pentagon*, color=blue, opacity=0.5, mark size=2] coordinates {(539.5,103.4)};
    \addlegendentry{$\bd^*$, 0.38 J }
    \addplot[only marks, mark=triangle*, color=blue, opacity=0.5, mark size=3] coordinates {(637.4,119.6)};
    \addlegendentry{$\bd^*$, 0.38 J }

    \addplot[dashed, black] coordinates {(100,20) (1000,20)};
    \addplot[dashed, black] coordinates {(100,200) (1000,200)};
    \addplot[dashed, black] coordinates {(100,20) (100,200)};
    \addplot[dashed, black] coordinates {(1000,20) (1000,200)};

    \node[anchor=north east, fill=white, draw, dashed] at (rel axis cs:0.935, 0.265) {Design space};
    \end{axis}
\end{tikzpicture}

%% file: main.bbl
\begin{thebibliography}{10}

\bibitem{bartlett2019overview}
{\sc Bartlett, J.~L., and Li, X.}
\newblock An overview of residual stresses in metal powder bed fusion.
\newblock {\em Additive Manufacturing 27\/} (2019), 131--149.

\bibitem{bhardwaj2019direct}
{\sc Bhardwaj, T., Shukla, M., Paul, C., and Bindra, K.}
\newblock Direct energy deposition-laser additive manufacturing of
  titanium-molybdenum alloy: Parametric studies, microstructure and mechanical
  properties.
\newblock {\em Journal of Alloys and Compounds 787\/} (2019), 1238--1248.

\bibitem{cansizoglu2008applications}
{\sc Cansizoglu, O., Harrysson, O.~L., West, H.~A., Cormier, D.~R., and Mahale,
  T.}
\newblock Applications of structural optimization in direct metal fabrication.
\newblock {\em Rapid Prototyping Journal 14}, 2 (2008), 114--122.

\bibitem{cao2021optimization}
{\sc Cao, L., Li, J., Hu, J., Liu, H., Wu, Y., and Zhou, Q.}
\newblock Optimization of surface roughness and dimensional accuracy in {LPBF}
  additive manufacturing.
\newblock {\em Optics \& Laser Technology 142\/} (2021), 107246.

\bibitem{chastand2018comparative}
{\sc Chastand, V., Quaegebeur, P., Maia, W., and Charkaluk, E.}
\newblock Comparative study of fatigue properties of {Ti-6Al-4V} specimens
  built by electron beam melting ({EBM}) and selective laser melting ({SLM}).
\newblock {\em Materials Characterization 143\/} (2018), 76--81.

\bibitem{chaudhuri2022certifiable}
{\sc Chaudhuri, A., Kramer, B., Norton, M., Royset, J.~O., and Willcox, K.}
\newblock Certifiable risk-based engineering design optimization.
\newblock {\em AIAA Journal 60}, 2 (2022), 551--565.

\bibitem{chia2022process}
{\sc Chia, H.~Y., Wu, J., Wang, X., and Yan, W.}
\newblock Process parameter optimization of metal additive manufacturing: A
  review and outlook.
\newblock {\em Journal of Materials Informatics 2}, 4 (2022).

\bibitem{constantine2015active}
{\sc Constantine, P.~G.}
\newblock {\em Active subspaces}.
\newblock Society for Industrial and Applied Mathematics, Philadelphia, PA,
  2015.

\bibitem{corbin2018effect}
{\sc Corbin, D.~J., Nassar, A.~R., Reutzel, E.~W., Beese, A.~M., and
  Michaleris, P.}
\newblock Effect of substrate thickness and preheating on the distortion of
  laser deposited {Ti--6Al--4V}.
\newblock {\em Journal of Manufacturing Science and Engineering 140}, 6 (2018),
  061009.

\bibitem{debroy2017building}
{\sc Debroy, T., Zhang, W., Turner, J., and Babu, S.~S.}
\newblock Building digital twins of {3D} printing machines.
\newblock {\em Scripta Materialia 135\/} (2017), 119--124.

\bibitem{elsayed2018optimization}
{\sc Elsayed, M., Ghazy, M., Youssef, Y., and Essa, K.}
\newblock Optimization of {SLM} process parameters for {Ti6Al4V} medical
  implants.
\newblock {\em Rapid Prototyping Journal 25}, 3 (2018), 433--447.

\bibitem{everton2016review}
{\sc Everton, S.~K., Hirsch, M., Stravroulakis, P., Leach, R.~K., and Clare,
  A.~T.}
\newblock Review of in-situ process monitoring and in-situ metrology for metal
  additive manufacturing.
\newblock {\em Materials \& Design 95\/} (2016), 431--445.

\bibitem{fu20143}
{\sc Fu, C., and Guo, Y.}
\newblock 3-dimensional finite element modeling of selective laser melting
  {Ti-6Al-4V} alloy.
\newblock In {\em Proceedings of the 2014 International Solid Freeform
  Fabrication Symposium\/} (2014), University of Texas at Austin,
  pp.~1129--1144.

\bibitem{gaikwad2020heterogeneous}
{\sc Gaikwad, A., Giera, B., Guss, G.~M., Forien, J.-B., Matthews, M.~J., and
  Rao, P.}
\newblock Heterogeneous sensing and scientific machine learning for quality
  assurance in laser powder bed fusion--a single-track study.
\newblock {\em Additive Manufacturing 36\/} (2020), 101659.

\bibitem{galarraga2016effects}
{\sc Galarraga, H., Lados, D.~A., Dehoff, R.~R., Kirka, M.~M., and Nandwana,
  P.}
\newblock Effects of the microstructure and porosity on properties of
  {Ti-6Al-4V} {ELI} alloy fabricated by electron beam melting ({EBM}).
\newblock {\em Additive Manufacturing 10\/} (2016), 47--57.

\bibitem{gibson2021additive}
{\sc Gibson, I., Rosen, D.~W., Stucker, B., Khorasani, M., Rosen, D., Stucker,
  B., and Khorasani, M.}
\newblock {\em Additive manufacturing technologies}, vol.~17.
\newblock Springer, 2021.

\bibitem{gong2015influence}
{\sc Gong, H., Rafi, K., Gu, H., Ram, G.~J., Starr, T., and Stucker, B.}
\newblock Influence of defects on mechanical properties of {Ti--6Al--4V}
  components produced by selective laser melting and electron beam melting.
\newblock {\em Materials \& Design 86\/} (2015), 545--554.

\bibitem{gong2014analysis}
{\sc Gong, H., Rafi, K., Gu, H., Starr, T., and Stucker, B.}
\newblock Analysis of defect generation in {Ti--6Al--4V} parts made using
  powder bed fusion additive manufacturing processes.
\newblock {\em Additive Manufacturing 1\/} (2014), 87--98.

\bibitem{cdi_astm_standards_111828}
{\sc International Organization for Standardization, American Society for
  Testing and Materials}.
\newblock {\em {ISO/ASTM52900} {A}dditive manufacturing - {G}eneral principles
  - {F}undamentals and vocabulary}, 2021.
\newblock Version F3177-21.

\bibitem{kalentics2017tailoring}
{\sc Kalentics, N., Boillat, E., Peyre, P., {\'C}iri{\'c}-Kosti{\'c}, S.,
  Bogojevi{\'c}, N., and Log{\'e}, R.~E.}
\newblock Tailoring residual stress profile of selective laser melted parts by
  laser shock peening.
\newblock {\em Additive Manufacturing 16\/} (2017), 90--97.

\bibitem{karlsson2013characterization}
{\sc Karlsson, J., Snis, A., Engqvist, H., and Lausmaa, J.}
\newblock Characterization and comparison of materials produced by electron
  beam melting ({EBM}) of two different {Ti--6Al--4V} powder fractions.
\newblock {\em Journal of Materials Processing Technology 213}, 12 (2013),
  2109--2118.

\bibitem{konda2017additive}
{\sc Konda~Gokuldoss, P., Kolla, S., and Eckert, J.}
\newblock Additive manufacturing processes: Selective laser melting, electron
  beam melting and binder jetting-selection guidelines.
\newblock {\em Materials 10}, 6 (2017), 672.

\bibitem{lee2021robust}
{\sc Lee, D., and Rahman, S.}
\newblock Robust design optimization under dependent random variables by a
  generalized polynomial chaos expansion.
\newblock {\em Structural and Multidisciplinary Optimization 63}, 5 (2021),
  2425--2457.

\bibitem{li2018residual}
{\sc Li, C., Liu, Z., Fang, X., and Guo, Y.}
\newblock Residual stress in metal additive manufacturing.
\newblock {\em Procedia {CIRP} 71\/} (2018), 348--353.

\bibitem{lu2021additive}
{\sc Lu, W., Zhai, W., Wang, J., Liu, X., Zhou, L., Ibrahim, A. M.~M., Li, X.,
  Lin, D., and Wang, Y.~M.}
\newblock Additive manufacturing of isotropic-grained, high-strength and
  high-ductility copper alloys.
\newblock {\em Additive Manufacturing 38\/} (2021), 101751.

\bibitem{lu2015study}
{\sc Lu, Y., Wu, S., Gan, Y., Huang, T., Yang, C., Junjie, L., and Lin, J.}
\newblock Study on the microstructure, mechanical property and residual stress
  of {SLM Inconel-718} alloy manufactured by differing island scanning
  strategy.
\newblock {\em Optics \& Laser Technology 75\/} (2015), 197--206.

\bibitem{megahed2016metal}
{\sc Megahed, M., Mindt, H.-W., N'Dri, N., Duan, H., and Desmaison, O.}
\newblock Metal additive-manufacturing process and residual stress modeling.
\newblock {\em Integrating Materials and Manufacturing Innovation 5}, 1 (2016),
  61--93.

\bibitem{meng2020process}
{\sc Meng, L., and Zhang, J.}
\newblock Process design of laser powder bed fusion of stainless steel using a
  gaussian process-based machine learning model.
\newblock {\em JOM 72}, 1 (2020), 420--428.

\bibitem{mur1996influence}
{\sc Mur, F.~G., Rodr{\'\i}guez, D., and Planell, J.}
\newblock Influence of tempering temperature and time on the
  $\alpha$'-{Ti-6Al-4V} martensite.
\newblock {\em Journal of Alloys and Compounds 234}, 2 (1996), 287--289.

\bibitem{murr2018additive}
{\sc Murr, L.~E.}
\newblock Additive manufacturing of biomedical devices: an overview.
\newblock {\em Materials Technology 33}, 1 (2018), 57--70.

\bibitem{oyesola2021optimization}
{\sc Oyesola, M., Mpofu, K., Mathe, N., Fatoba, S., Hoosain, S., and Daniyan,
  I.}
\newblock Optimization of selective laser melting process parameters for
  surface quality performance of the fabricated {Ti6Al4V}.
\newblock {\em The International Journal of Advanced Manufacturing Technology
  114\/} (2021), 1585--1599.

\bibitem{rockafellar2010buffered}
{\sc Rockafellar, R.~T., and Royset, J.~O.}
\newblock On buffered failure probability in design and optimization of
  structures.
\newblock {\em Reliability Engineering \& System Safety 95}, 5 (2010),
  499--510.

\bibitem{rockafellar2002conditional}
{\sc Rockafellar, R.~T., and Uryasev, S.}
\newblock Conditional value-at-risk for general loss distributions.
\newblock {\em Journal of Banking \& Finance 26}, 7 (2002), 1443--1471.

\bibitem{seoane2021semi}
{\sc Seoane-Via{\~n}o, I., Januskaite, P., Alvarez-Lorenzo, C., Basit, A.~W.,
  and Goyanes, A.}
\newblock Semi-solid extrusion {3D} printing in drug delivery and biomedicine:
  Personalised solutions for healthcare challenges.
\newblock {\em Journal of Controlled Release 332\/} (2021), 367--389.

\bibitem{shapiro2016additive}
{\sc Shapiro, A.~A., Borgonia, J., Chen, Q., Dillon, R., McEnerney, B.,
  Polit-Casillas, R., and Soloway, L.}
\newblock Additive manufacturing for aerospace flight applications.
\newblock {\em Journal of Spacecraft and Rockets\/} (2016), 952--959.

\bibitem{shi2017parameter}
{\sc Shi, X., Ma, S., Liu, C., and Wu, Q.}
\newblock Parameter optimization for {Ti-47Al-2Cr-2Nb} in selective laser
  melting based on geometric characteristics of single scan tracks.
\newblock {\em Optics \& Laser Technology 90\/} (2017), 71--79.

\bibitem{shiomi2004residual}
{\sc Shiomi, M., Osakada, K., Nakamura, K., Yamashita, T., and Abe, F.}
\newblock Residual stress within metallic model made by selective laser melting
  process.
\newblock {\em CIRP Annals 53}, 1 (2004), 195--198.

\bibitem{spears2016process}
{\sc Spears, T.~G., and Gold, S.~A.}
\newblock In-process sensing in selective laser melting ({SLM}) additive
  manufacturing.
\newblock {\em Integrating Materials and Manufacturing Innovation 5}, 1 (2016),
  16--40.

\bibitem{sun2013parametric}
{\sc Sun, J., Yang, Y., and Wang, D.}
\newblock Parametric optimization of selective laser melting for forming
  {Ti6Al4V} samples by {Taguchi} method.
\newblock {\em Optics \& Laser Technology 49\/} (2013), 118--124.

\bibitem{vasco2021additive}
{\sc Vasco, J.~C.}
\newblock Additive manufacturing for the automotive industry.
\newblock In {\em Additive manufacturing}. Elsevier, 2021, pp.~505--530.

\bibitem{vastola2016controlling}
{\sc Vastola, G., Zhang, G., Pei, Q., and Zhang, Y.-W.}
\newblock Controlling of residual stress in additive manufacturing of {Ti6Al4V}
  by finite element modeling.
\newblock {\em Additive Manufacturing 12\/} (2016), 231--239.

\bibitem{vilaro2011fabricated}
{\sc Vilaro, T., Colin, C., and Bartout, J.-D.}
\newblock As-fabricated and heat-treated microstructures of the {Ti-6Al-4V}
  alloy processed by selective laser melting.
\newblock {\em Metallurgical and materials transactions A 42}, 10 (2011),
  3190--3199.

\bibitem{2020SciPy-NMeth}
{\sc Virtanen, P., Gommers, R., Oliphant, T.~E., Haberland, M., Reddy, T.,
  Cournapeau, D., Burovski, E., Peterson, P., Weckesser, W., Bright, J., {van
  der Walt}, S.~J., Brett, M., Wilson, J., Millman, K.~J., Mayorov, N., Nelson,
  A. R.~J., Jones, E., Kern, R., Larson, E., Carey, C.~J., Polat, {\.I}., Feng,
  Y., Moore, E.~W., {VanderPlas}, J., Laxalde, D., Perktold, J., Cimrman, R.,
  Henriksen, I., Quintero, E.~A., Harris, C.~R., Archibald, A.~M., Ribeiro,
  A.~H., Pedregosa, F., {van Mulbregt}, P., and {SciPy 1.0 Contributors}}.
\newblock {{SciPy} 1.0: Fundamental Algorithms for Scientific Computing in
  Python}.
\newblock {\em Nature Methods 17\/} (2020), 261--272.

\bibitem{wang2024uncertainty}
{\sc Wang, H., Li, B., Lei, L., and Xuan, F.}
\newblock Uncertainty-aware fatigue-life prediction of additively manufactured
  {Hastelloy X} superalloy using a physics-informed probabilistic neural
  network.
\newblock {\em Reliability Engineering \& System Safety 243\/} (2024), 109852.

\bibitem{wang2019data}
{\sc Wang, Z., Liu, P., Xiao, Y., Cui, X., Hu, Z., and Chen, L.}
\newblock A data-driven approach for process optimization of metallic additive
  manufacturing under uncertainty.
\newblock {\em Journal of Manufacturing Science and Engineering 141}, 8 (2019),
  081004.

\bibitem{wei2018heterogeneous}
{\sc Wei, D., Koizumi, Y., Chiba, A., Ueki, K., Ueda, K., Narushima, T.,
  Tsutsumi, Y., and Hanawa, T.}
\newblock Heterogeneous microstructures and corrosion resistance of biomedical
  {Co-Cr-Mo} alloy fabricated by electron beam melting ({EBM}).
\newblock {\em Additive Manufacturing 24\/} (2018), 103--114.

\bibitem{welsch1993materials}
{\sc Welsch, G., Boyer, R., and Collings, E.}
\newblock {\em Materials properties handbook: titanium alloys}.
\newblock ASM international, 1993.

\bibitem{wu2014experimental}
{\sc Wu, A.~S., Brown, D.~W., Kumar, M., Gallegos, G.~F., and King, W.~E.}
\newblock An experimental investigation into additive manufacturing-induced
  residual stresses in {316L} stainless steel.
\newblock {\em Metallurgical and Materials Transactions A 45\/} (2014),
  6260--6270.

\bibitem{yadroitsev2013energy}
{\sc Yadroitsev, I., Krakhmalev, P., Yadroitsava, I., Johansson, S., and
  Smurov, I.}
\newblock Energy input effect on morphology and microstructure of selective
  laser melting single track from metallic powder.
\newblock {\em Journal of Materials Processing Technology 213}, 4 (2013),
  606--613.

\bibitem{yan2017multi}
{\sc Yan, W., Ge, W., Qian, Y., Lin, S., Zhou, B., Liu, W.~K., Lin, F., and
  Wagner, G.~J.}
\newblock Multi-physics modeling of single/multiple-track defect mechanisms in
  electron beam selective melting.
\newblock {\em Acta Materialia 134\/} (2017), 324--333.

\bibitem{zhao2015numerical}
{\sc Zhao, X., Promoppatum, P., and Yao, S.-C.}
\newblock Numerical modeling of non-linear thermal stress in direct metal laser
  sintering process of titanium alloy products.
\newblock In {\em Proceedings of the First Thermal and Fluids Engineering
  Summer Conference\/} (2015), pp.~9--12.

\bibitem{zinoviev2016evolution}
{\sc Zinoviev, A., Zinovieva, O., Ploshikhin, V., Romanova, V., and Balokhonov,
  R.}
\newblock Evolution of grain structure during laser additive manufacturing.
  simulation by a cellular automata method.
\newblock {\em Materials \& Design 106\/} (2016), 321--329.

\bibitem{zou2006direct}
{\sc Zou, T., and Mahadevan, S.}
\newblock A direct decoupling approach for efficient reliability-based design
  optimization.
\newblock {\em Structural and Multidisciplinary Optimization 31\/} (2006),
  190--200.

\end{thebibliography}
